\documentclass[11pt]{amsart}
\usepackage{amsmath,amsfonts,amssymb,url,graphicx}

\DeclareFontFamily{OT1}{rsfs}{}
\DeclareFontShape{OT1}{rsfs}{n}{it}{<-> rsfs10}{}
\DeclareMathAlphabet{\mathscr}{OT1}{rsfs}{n}{it}
\DeclareMathOperator{\sgn}{sgn}

\DeclareMathOperator{\mo}{\,mod}

\newtheorem{prop}{Proposition}[section]
\newtheorem{thm}[prop]{Theorem}
\newtheorem*{main}{Main Theorem}

\newtheorem*{defn*}{Definition}

\numberwithin{equation}{section}

\title{The ternary Goldbach problem}
\author{Harald Andr\'es Helfgott}
\date{}
\address{Harald Helfgott, 
\'Ecole Normale Sup\'erieure, D\'epartement de Math\'ematiques, 45 rue d'Ulm, F-75230 Paris, France}
\email{harald.helfgott@ens.fr}
\begin{document}
\begin{abstract}
The ternary Goldbach conjecture, or three-primes problem, states
that every odd number $n$ greater than $5$ can be written as the sum
of three primes. The conjecture, posed in 1742, remained unsolved until
now, in spite of great progress in the twentieth century. In 2013
-- following a line of research pioneered and developed by Hardy, Littlewood and
Vinogradov, among others -- the author proved the conjecture.

In this, as in many other additive problems, what is at issue is really
the proper usage of the limited information we possess on the
distribution of prime numbers. The problem serves as a test and
whetting-stone for techniques in analysis and number theory -- and
also as an incentive to think about the relations between existing
techniques with greater clarity.

We will go over the main ideas of the proof. The basic approach is based
on the circle method, the large sieve and exponential sums. For the purposes
of this overview, we will not need to work with explicit constants; however,
we will discuss what makes certain strategies and procedures
not just effective, but efficient, in the sense of leading to good constants.
Still, our focus will be on qualitative improvements.
\end{abstract}

\maketitle

The question we will discuss, or one similar to it, seems to have been
first posed by Descartes, in a manuscript published only centuries after
his death \cite[p. 298]{zbMATH02639585}. Descartes states:
``Sed \& omnis numerus par fit ex uno vel duobus vel tribus primis''
(``But also every even number is made out of one, two or three prime numbers.''\footnote{Thanks are due to J. Brandes and R. Vaughan for a discussion on a possible
 ambiguity in the Latin wording. Descartes' statement is 
mentioned (with a translation much like the one given here) in Dickson's
{\em History} \cite[Ch. XVIII]{MR0245499}.}.)
This statement comes in the middle of a discussion of sums of 
polygonal numbers, such as the squares.

Statements on sums of primes and sums of values of polynomials 
(polygonal numbers, powers $n^k$, etc.) 
have since shown themselves to be much more
than mere curiosities -- and not just because they are often very difficult
to prove. Whereas the study of sums of powers can rely on
 their algebraic structure, the study of sums of primes leads to the 
realization that, from several perspectives, the set of primes behaves much 
like the set of integers -- and that this is truly hard to prove. 

If, instead of the primes, we had a random set of odd integers $S$ whose density
-- an intuitive concept that can be made precise -- equaled that of the 
primes, then we would expect to be able to write every odd number as a sum
of three elements of $S$, and every even number as the sum of two elements
of $S$. We would have to check by hand whether this is true for small odd
and even numbers, but it is relatively easy to show that,
after a long enough check, it would be very unlikely
that there would be any exceptions left among the infinitely many cases left
to check.

The question, then, is in what sense we need the primes to be like a random
set of integers;
in other words, we need to know what we can prove about the regularities of
the distribution of the primes. This is one of the main questions of analytic
number theory; progress on it has been very slow and difficult. Thus, the real
question is how to use well the limited information we do have on the 
distribution of the primes.


\section{History and new developments}

The history of the conjecture starts properly with Euler and 
his close friend, Christian Goldbach, both of whom lived and worked in Russia
at the time of their correspondence -- about a century after Descartes'
isolated statement. 
Goldbach, a man of many interests,
is usually classed as a serious amateur; he seems to have awakened Euler's 
passion for number theory, which would lead to the 
beginning of the modern era of the subject \cite[Ch. 3, \S IV]{MR734177}.
In a letter dated June 7, 1742 -- written partly in German, partly in Latin --
Goldbach made a conjectural statement on prime numbers, and Euler rapidly
reduced it to the following conjecture, which, he said, Goldbach had already
posed to him: every positive integer can be written as the sum of at most
three prime numbers.

We would now say ``every integer greater than $1$'', since we no long
consider $1$ to be a prime number. Moreover, the conjecture is nowadays split 
into two: 
\begin{itemize}
\item the {\em weak}, or ternary, Goldbach conjecture states that every odd integer greater than $5$ can be written as the sum of three primes; 
\item the {\em strong}, or binary, Goldbach conjecture states that every even integer greater than $2$ can be written as the sum of two primes.
\end{itemize}
As their names indicate, the strong conjecture implies the weak one 
(easily: subtract $3$ from your odd number $n$, then express $n-3$ as the sum 
of two primes).

The strong conjecture remains out of reach. A short while ago -- the first
complete version appeared on May 13, 2013 -- the present author
 proved the weak Goldbach conjecture.

\begin{main}
Every odd integer
greater than $5$ can be written as the sum of three primes.
\end{main}

The proof is contained in the preprints \cite{Helf}, \cite{HelfMaj}, 
\cite{HelfTern}. 
It builds on the great progress towards the conjecture
made in the early 20th century by Hardy, Littlewood and Vinogradov.
In 1937, Vinogradov proved \cite{zbMATH02522879} that the conjecture is true for all odd numbers $n$ 
larger than some constant $C$. 
(Hardy and Littlewood had shown the same under the assumption of the 
Generalized Riemann Hypothesis, which we shall have the chance to discuss
later.) 

It is clear that a computation can verify the conjecture only for $n\leq c$,
$c$ a constant: computations have to be finite. What can make a result
coming from analytic number theory be valid only for $n\geq C$?

An analytic proof, generally speaking, gives us more than just existence. In this kind of problem, it gives us more than the possibility of doing
something (here, writing an integer $n$ as the sum of three
primes). It gives us a rigorous estimate for the number of ways in 
which this {\em something}
is possible; that is, it shows us that this number of ways equals
\begin{equation}\label{eq:huta}
\text{main term} + \text{error term},\end{equation}
where the main term is a precise quantity $f(n)$, and the error term is  
something whose absolute value is at most another precise quantity $g(n)$.
If $f(n)>g(n)$, then (\ref{eq:huta}) is non-zero, i.e., we will have
shown that the existence of a way to write our number as the sum of
three primes.

(Since what we truly care about is existence, we are free to weigh different
ways of writing $n$ as the sum of three primes however we wish -- that is,
we can decide that some primes ``count'' twice or thrice as much as others,
and that some do not count at all.)

Typically, after much work, we succeed in obtaining (\ref{eq:huta})
with $f(n)$ and $g(n)$ such that $f(n)>g(n)$ asymptotically, that is,
for $n$ large enough. To give a highly simplified example:
if, say, $f(n) = n^2$ and $g(n) = 100 n^{3/2}$, then $f(n)>g(n)$ for
$n>C$, where $C=10^4$, and so the number of ways
(\ref{eq:huta}) is positive for $n>C$.

We want a moderate value of $C$, that is, a $C$ small enough that all cases
$n\leq C$ can be checked computationally. To ensure this, we must make
the error term bound $g(n)$ as small as possible. This is our main
task. A secondary (and sometimes neglected)
 possibility is to rig the weights so as to make the main term
$f(n)$ larger in comparison to $g(n)$; this can generally be done only up to a
certain point, but is nonetheless very helpful.

As we said, the first unconditional proof that odd numbers $n\geq C$
can be written as the sum of three primes is due to Vinogradov.
Analytic bounds fall into several categories, or stages; quite often, successive
versions of the same theorem will go through successive stages.
\begin{enumerate}
\item An {\em ineffective} result shows that a statement is true for some
constant $C$, but gives no way to determine what the constant $C$ might be.
Vinogradov's first proof of his theorem (in \cite{zbMATH02522879}) is like this:
it shows that there exists a constant $C$ such that
 every odd number $n>C$ is the sum of three primes, yet give us no hope
of finding out what the constant $C$ might be.\footnote{Here, as is often the 
case in ineffective results in analytic number theory, the underlying 
issue is that of {\em Siegel zeros}, which are believed not to exist, but
have not been shown not to; the strongest bounds on (i.e., against) such zeros
are ineffective, and so are all of the many results using such estimates.}
Many proofs of Vinogradov's result in textbooks are also of this type.
\item An {\em effective}, but not explicit, result shows that a statement
is true for some unspecified constant $C$ in a way that makes it clear that
a constant $C$ could in principle be determined following and reworking
the proof with great care. Vinogradov's later proof 
(\cite{zbMATH03063033}, translated in \cite{MR0062183}) is of this nature.
As Chudakov \cite[\S IV.2]{MR0031961} pointed out, the improvement
on \cite{zbMATH02522879} given
by Mardzhanishvili \cite{Mardzh} already had the effect of making
the result effective.\footnote{The proof in \cite{Mardzh} 
combined the bounds in \cite{zbMATH02522879} with a more careful accounting of
the effect of the single possible Siegel zero within range.}
\item An {\em explicit} result gives a value of $C$. 
According to \cite[p. 201]{MR0031961}, the first explicit
  version of Vinogradov's result was given by Borozdkin in his unpublished
doctoral dissertation, written under the direction of Vinogradov 
(1939): $C = \exp(\exp(\exp(41.96)))$. Such a result is, by definition,
also effective. Borodzkin later
  \cite{Borodzkin} 
gave the value $C = e^{e^{16.038}}$, though he does not seem to have published
the proof. The best -- that is, smallest -- value of $C$ known before the 
present work was that of Liu and 
Wang \cite{MR1932763}: $C = 2\cdot 10^{1346}$.
\item What we may call an {\em efficient} proof gives a reasonable value
for $C$ -- in our case, a value small enough that checking all cases up to
$C$ is feasible.
\end{enumerate}

How far were we from an efficient proof? That is, what sort of computation 
could ever be feasible? The number of picoseconds since the beginning
of the universe is less than $10^{30}$, whereas the number of protons
in the observable universe is currently estimated at $\sim 10^{80}$
\cite{EB}. This means that even a parallel computer the size of the universe
could never perform a computation requiring $10^{110}$ steps, even if it ran
for the age of the universe. Thus, $C = 2\cdot 10^{1346}$ is too large.

I gave a proof with $C=10^{29}$ in May 2013. Since D. Platt and I had verified
the conjecture for all odd numbers up to $n\leq 8.8\cdot 10^{30}$
by computer \cite{HelPlat}, this established the conjecture for all
odd numbers $n$.

(In December 2013, $C$ was reduced
 to $10^{27}$ \cite{HelfTern}. The verification
of the ternary Goldbach conjecture up to $n\leq 10^{27}$ can be done
in a home computer over a weekend. All must be said: this uses 
the verification of the binary Goldbach conjecture for $n\leq 4\cdot 10^{18}$
\cite{OSHP}, which itself required computational resources far outside the 
home-computing range. 
Checking the conjecture up to $n\leq 10^{27}$ was not even 
 the main computational task that needed to be accomplished to establish
the Main Theorem -- that task was
the finite verification of zeros of $L$-functions in 
\cite{Plattfresh}, 
a general-purpose computation that should be useful elsewhere. 
We will discuss the procedure at the end of the
article.)

What was the strategy of \cite{Helf}, \cite{HelfMaj}, 
and \cite{HelfTern}? The basic framework is the one pioneered by
Hardy and Littlewood for a variety of problems -- namely, the
{\em circle method}, which, as we shall see, is an application of
Fourier analysis over $\mathbb{Z}$. 
(There are other, later routes to Vinogradov's result;
see \cite{MR834356}, \cite{MR1670069} 
and especially the recent work \cite{MR3165421}, which avoids using anything
about zeros of $L$-functions inside the critical strip.) 
Vinogradov's proof, 
like much of the later work on the subject, was based on a detailed analysis of
exponential sums, i.e., Fourier transforms over $\mathbb{Z}$. So is
the proof that we will sketch.

At the same time, the distance between $2\cdot 10^{1346}$ and $10^{27}$
is such that we cannot hope to get to $10^{27}$ (or any other
reasonable constant) by fine-tuning previous work.
Rather, we must work from scratch, using the basic outline
in Vinogradov's original proof and other, initially unrelated, developments
in analysis and number theory (notably, the large sieve). Merely improving
 constants will not do; rather, we must do qualitatively better than
previous work (by non-constant factors) if we are to have any chance to
succeed. It is on these qualitative improvements that we
will focus.

\begin{center}
* * *
\end{center}

It is only fair to review some of the progress made between Vinogradov's
time and ours. Here we will focus on results; later, 
we will discuss some of the progress made in the techniques of proof.
For a fuller account up to 1978, see R. Vaughan's
ICM lecture notes on the ternary Goldbach problem \cite{MR562631}.

In 1933, Schnirelmann proved 
\cite{MR1512821} that every integer $n>1$ can be written as the sum
of at most $K$ primes for some unspecified constant $K$. (This pioneering work
is now considered to be
part of the early history of additive combinatorics.) 
In 1969, Klimov gave an explicit value for $K$ (namely, 
$K = 6\cdot 10^9$); he later improved the constant to
$K = 115$ (with G. Z. Piltay and T. A. Sheptickaja) and $K = 55$. 
Later, there were results by Vaughan \cite{MR0437478} ($K=27$), 
Deshouillers
\cite{MR0466050} ($K=26$) and Riesel-Vaughan \cite{MR706639} ($K=19$).

Ramar\'e showed in 1995 that every even number $n>1$ can be written
as the sum of at most $6$ primes \cite{MR1375315}. In 2012, Tao proved
\cite{Tao} that every odd number $n>1$ is the sum of at most $5$ primes.

There have been other avenues of attack towards the strong conjecture.
Using ideas close to those of Vinogradov's,
Chudakov \cite{zbMATH03028355}, 
\cite{Chudtoo}, Estermann  \cite{MR1576891} and
van der Corput \cite{zbMATH02522863} 
 proved (independently from each other) that almost every even number
(meaning: all elements of a subset of density $1$ in the even numbers)
can be written as the sum of two primes. In 1973, J.-R. Chen
showed \cite{MR0434997} that every even number $n$ larger than a constant
$C$ can be written as the sum of a prime number and the product of
at most two primes ($n = p_1 + p_2$ or 
$n = p_1 + p_2 p_3$). Incidentally, J.-R. Chen himself, together with
T.-Z. Wang, was responsible for the best bounds on $C$ (for ternary
 Goldbach) before Lui and Wang: 
$C = \exp(\exp(11.503)) < 4\cdot 10^{43000}$
\cite{MR1046491} and
$C = \exp(\exp(9.715)) < 6\cdot
10^{7193}$ \cite{MR1411958}.

Matters are different if one assumes the Generalized Riemann
Hypothesis (GRH). A careful analysis \cite{MR1715106}
of Hardy and Littlewood's work \cite{MR1555183} gives that every 
odd number $n\geq 1.24\cdot 10^{50}$ is the sum of three primes if GRH
is true.
According to \cite{MR1715106}, the same statement with
 $n\geq 10^{32}$ was proven in the unpublished doctoral dissertation
of B. Lucke, a student of E. Landau's, in 1926. Zinoviev \cite{MR1462848}
improved this to $n\geq 10^{20}$. A computer check (\cite{MR1469323}; see also
\cite{MR1451327}) showed that the conjecture is true for $n<10^{20}$,
thus completing the proof of the ternary Goldbach conjecture 
under the assumption of GRH. What was open until now was, of course,
the problem of giving an unconditional proof.
 
{\bf Acknowledgments.} 
Parts of the present article are based on a previous expository note by the 
author. The first version of the note
appeared online, in English, in an informal venue
\cite{Helblog}; later versions were
published in Spanish (\cite{Gaceta}, translated by M. A. Morales and the author,
and revised with the help of J. Cilleruelo and M. Helfgott)
 and French (\cite{Gazette}, translated by M. Bilu and revised by the author). 
Many individuals and organizations should be thanked for their generous help towards the work
summarized here; an attempt at a full list can be found in the acknowledgments
sections of  \cite{Helf}, \cite{HelfMaj}, \cite{HelfTern}. 
Thanks are also due to J. Brandes, K. Gong, R. Heath-Brown,
Z. Silagadze, R. Vaughan and T. Wooley, for help with 
historical questions.

\section{The circle method: Fourier analysis on $\mathbb{Z}$}

It is common for a first course on Fourier analysis to focus
on functions over the reals satisfying $f(x)=f(x+1)$, or, what is the same,
functions $f:\mathbb{R}/\mathbb{Z} \to \mathbb{C}$. Such a function
(unless it is fairly pathological) has a Fourier series converging to it;
this is just the same as saying that $f$ has a Fourier transform
$\widehat{f}:\mathbb{Z}\to \mathbb{C}$ defined by 
$\widehat{f}(n) = \int_{\mathbb{R}/\mathbb{Z}} f(\alpha) e(-\alpha n) d\alpha$ 
and satisfying $f(\alpha) = \sum_{n\in \mathbb{Z}} \widehat{f}(n) e(\alpha
n) d\alpha$
({\em Fourier inversion theorem}). 

In number theory, we are especially interested in functions $f:\mathbb{Z}
\to \mathbb{C}$. Then things are exactly the other way around: provided
that $f$ decays reasonably fast as $n\to \pm \infty$ (or becomes $0$ for $n$ large enough), 
$f$ has a Fourier transform $\widehat{f}:\mathbb{R}/\mathbb{Z} \to
\mathbb{C}$ defined by 
$\widehat{f}(\alpha) = \sum_n f(n) e(-\alpha n)$
and satisfying $f(n) = \int_{\mathbb{R}/\mathbb{Z}} \widehat{f}(\alpha) e(\alpha n)$.
(Highbrow talk: we already knew that $\mathbb{Z}$ is the Fourier dual of
$\mathbb{R}/\mathbb{Z}$, and so, of course, 
$\mathbb{R}/\mathbb{Z}$ is the Fourier dual of $\mathbb{Z}$.) 
``Exponential sums'' (or ``trigonometrical sums'', as in the title of
\cite{MR0062183})
are sums of the form $\sum_n f(\alpha) e(-\alpha n)$; the ``circle''
in ``circle method'' is just a name for $\mathbb{R}/\mathbb{Z}$.

The study of the Fourier transform $\widehat{f}$ 
is relevant to additive problems in number theory, i.e., questions
on the number of ways of writing $n$ as a sum of $k$ integers of a particular
form. Why? One answer could be that $\widehat{f}$ gives us information
about the ``randomness'' of $f$; if $f$ were the characteristic function
of a random set, then $\widehat{f}(\alpha)$ would be very small outside a sharp
peak at $\alpha=0$.

We can also give a more concrete and immediate answer. Recall that, in 
general, the Fourier transform of a convolution
equals the product of the transforms; over $\mathbb{Z}$, this means that
for the additive convolution
\[(f\ast g)(n) = \mathop{\sum_{m_1, m_2\in \mathbb{Z}}}_{m_1+m_2 = n} f(m_1) g(m_2),\]
the Fourier transform satisfies the simple rule
\[\widehat{f\ast g}(\alpha) = \widehat{f}(\alpha) \cdot \widehat{g}(\alpha).\]

We can see right away from this that $(f\ast g)(n)$ can be non-zero only if $n$ 
can be written as $n=m_1+m_2$ for some $m_1$, $m_2$ such that $f(m_1)$ and 
$g(m_2)$ are non-zero. Similarly, $(f\ast g\ast h)(n)$ can be non-zero only if 
$n$ can be written as $n=m_1+m_2+m_3$ for some $m_1$, $m_2$, $m_3$ such that 
$f(m_1)$, $f_2(m_2)$ and $f_3(m_3)$ are all non-zero. This suggests that, to 
study the ternary Goldbach problem, we define $f_1, f_2, f_3: \mathbb{Z}\to
\mathbb{C}$ so that they take non-zero values only at the primes.

Hardy and Littlewood defined $f_1(n) = f_2(n) = f_3(n) = 0$ for $n$ non-prime 
(and also for $n\leq 0$), and $f_1(n) = f_2(n) = f_3(n) = (\log n) e^{-n/x}$ 
for $n$ prime (where $x$ is a parameter to be fixed later). 
Here the factor $e^{-n/x}$ is there to provide ``fast decay'', so that 
everything converges; as we will see later, Hardy and Littlewood's choice of 
$e^{-n/x}$ (rather than some other function of fast decay) is actually very 
clever, though not quite best-possible. The term $\log n$ is there for technical reasons -- in essence, it makes sense to put it there because a random
integer around $n$ has a chance of about $1/(\log n)$ of being prime.

We can see that $(f_1\ast f_2\ast f_3)(n) \ne 0$ if and only if $n$ can be written as  the sum of three primes. Our task is then to show that 
$(f_1\ast f_2\ast f_3)(n)$ 
(i.e., $(f\ast f\ast f)(n)$) is non-zero for every $n$ larger than a constant
$C\sim 10^{27}$. 
Since the transform of a convolution equals a product of 
transforms,
\begin{equation}\label{eq:vanes}
(f_1\ast f_2\ast f_3)(n) 
= \int_{\mathbb{R}/\mathbb{Z}} 
\widehat{f_1\ast f_2\ast f_3}(\alpha) e(\alpha n) d\alpha = 
\int_{\mathbb{R}/\mathbb{Z}} (\widehat{f_1} \widehat{f_2} \widehat{f_3})(\alpha) e(\alpha n) d\alpha.\end{equation}
 Our task is thus to show that the integral 
$\int_{\mathbb{R}/\mathbb{Z}} (\widehat{f_1} \widehat{f_2} \widehat{f_3})(\alpha) e(\alpha n) d\alpha$
is non-zero.

As it happens, $\widehat{f}(\alpha)$ is particularly large when $\alpha$ is close to a rational with small denominator.
Moreover, for such $\alpha$, it turns out we can actually give rather
precise estimates for $\widehat{f}(\alpha)$. Define $\mathfrak{M}$
(called the set of {\em major arcs}) to be a union of narrow arcs around the rationals with small denominator:
\[\mathfrak{M} = \bigcup_{q\leq r} \mathop{\bigcup_{a \mo q}}_{(a,q)=1}
\left(\frac{a}{q} - \frac{1}{q Q} , \frac{a}{q} + \frac{1}{q Q}\right),\]
where $Q$ is a constant times $x/r$, and $r$ will be set later.
We can write
\begin{equation}\label{eq:willow}
\int_{\mathbb{R}/\mathbb{Z}} (\widehat{f_1} \widehat{f_2} \widehat{f_3})(\alpha) e(\alpha n) d\alpha = 
\int_{\mathfrak{M}} (\widehat{f_1} \widehat{f_2} \widehat{f_3})(\alpha) e(\alpha n) d\alpha +
\int_{\mathfrak{m}} (\widehat{f_1} \widehat{f_2} \widehat{f_3})(\alpha) e(\alpha n) d\alpha,\end{equation}
where $\mathfrak{m}$ is the complement 
$(\mathbb{R}/\mathbb{Z})\setminus \mathfrak{M}$ (called {\em minor arcs}).

Now, we simply do not know how to give precise estimates for 
$\widehat{f}(\alpha)$ when $\alpha$ is in $\mathfrak{m}$. However, as
Vinogradov realized, one can give reasonable upper bounds on 
$|\widehat{f}(\alpha)|$ for $\alpha\in \mathfrak{m}$. This suggests
the following strategy: show that
\begin{equation}\label{eq:wtree}
\int_{\mathfrak{m}} |\widehat{f_1}(\alpha)| |\widehat{f_2}(\alpha)| 
|\widehat{f_3}(\alpha)| d\alpha 
<
\int_{\mathfrak{M}} \widehat{f_1}(\alpha) \widehat{f_2}(\alpha)
 \widehat{f_3}(\alpha) e(\alpha n) d\alpha.
\end{equation}
By (\ref{eq:vanes}) and (\ref{eq:willow}), this will imply immediately that
$(f_1\ast f_2\ast f_3)(n)  > 0$, and so we will be done.

The name of {\em circle method} is given to the study of additive problems
by means of Fourier analysis over $\mathbb{Z}$, and, in particular,
to the use of a subdivision of the circle $\mathbb{R}/\mathbb{Z}$ into
major and minor arcs to estimate the integral of a Fourier transform.
There was a ``circle'' already in Hardy and Ramanujan's work \cite{MR2280879},
but 
the subdivision into major and minor arcs is due to Hardy and Littlewood,
who also applied their method to a wide variety
of additive problems. (Hence ``the Hardy-Littlewood method'' as an alternative
name for the circle method.)
Before working on the ternary Goldbach conjecture, Hardy and Littlewood
also studied the question of whether every $n>C$ can be
written as the sum of $k$th powers, for instance. Vinogradov then 
showed how to do without contour integrals and worked with 
finite exponential sums, i.e., $f_i$ compactly supported. 
From today's perspective, it is clear that there are applications
(such as ours) in which it can be more important for $f_i$ to be
smooth than compactly supported; still, Vinogradov's simplifications
were a key incentive to further developments.

An important note: in the case of the binary Goldbach conjecture, the method
fails at (\ref{eq:wtree}), and not before; if our understanding of
the actual value of $\widehat{f_i}(\alpha)$ is at all correct, it is simply
not true in general that
\[\int_{\mathfrak{m}} |\widehat{f_1}(\alpha)| |\widehat{f_2}(\alpha)| d\alpha
<
\int_{\mathfrak{M}} \widehat{f_1}(\alpha) 
\widehat{f_2}(\alpha) e(\alpha n) d\alpha.\]
Let us see why this is not surprising. 
Set $f_1=f_2=f_3=f$ for simplicity, so that we have the integral of the square
$(\widehat{f}(\alpha))^2$
for the binary problem, and the integral of the cube 
$(\widehat{f}(\alpha))^3$ for the ternary problem.
Squaring, like cubing, amplifies the peaks of $\widehat{f}(\alpha)$, which are at
the rationals of small denominator and their immediate neighborhoods
(the major arcs); however, cubing amplifies the peaks much more than squaring.
This is why, even though the arcs making up $\mathfrak{M}$ are very narrow,
$\int_{\mathfrak{M}} f(\alpha)^3 e(\alpha n) d\alpha$ is larger than
$\int_{\mathfrak{m}} |f(\alpha)|^3 d\alpha$; that explains the name
{\em major arcs} -- they are not large, but they
give the major part of the contribution. In contrast, squaring
amplifies the peaks less, and this is why the absolute value of
$\int_{\mathfrak{M}} f(\alpha)^2 e(\alpha n) d\alpha$ is in general smaller than
$\int_{\mathfrak{m}} |f(\alpha)|^2 d\alpha$. As nobody knows how to prove
a precise estimate (and, in particular, lower bounds) on
$f(\alpha)$ for $\alpha \in \mathfrak{m}$, the binary Goldbach conjecture is
still very much out of reach.

To prove the ternary Goldbach conjecture, it is enough to estimate both sides
of (\ref{eq:wtree}) for carefully chosen $f_1$, $f_2$, $f_3$, and compare them.
This is our task from now on.

\section{The major arcs $\mathfrak{M}$}\label{sec:rubosa}
\subsection{What do we really know about $L$-functions and their zeros?}
Before we start, let us give 
a very brief review of basic analytic number theory
(in the sense of, say, \cite{MR0217022}). A {\em Dirichlet character}
$\chi:\mathbb{Z}\to \mathbb{C}$ of modulus $q$ is a character of 
$(\mathbb{Z}/q \mathbb{Z})^*$ lifted to $\mathbb{Z}$. (In other words,
$\chi(n) = \chi(n+q)$, $\chi(a b) = \chi(a) \chi(b)$ for all $a$, $b$
and $\chi(n)=0$ for $(n,q)\ne 1$.) A {\em Dirichlet $L$-series}
is defined by 
\[L(s,\chi) = \sum_{n=1}^\infty \chi(n) n^{-s}\]
for $\Re(s)>1$, and by analytic continuation for $\Re(s)\leq 1$. (The Riemann
zeta function $\zeta(s)$ is the $L$-function for the trivial character, i.e.,
the character $\chi$ such that $\chi(n)=1$ for all $n$.) Taking logarithms
and then derivatives, we see that
\begin{equation}\label{eq:koj}
- \frac{L'(s,\chi)}{L(s,\chi)} = \sum_{n=1}^\infty \Lambda(n) n^{-s},
\end{equation}
where $\Lambda$ is the {\em von Mangoldt function} ($\Lambda(n) = \log p$ if $n$ is some prime power $p^\alpha$, $\alpha\geq
1$, and $\Lambda(n)=0$ otherwise).

Dirichlet introduced his characters
and $L$-series so as to study primes in arithmetic progressions.
In general, and after some work,
(\ref{eq:koj}) allows us to restate many sums over the primes (such as our 
Fourier
transforms $\widehat{f}(\alpha)$) as sums over the zeros of $L(s,\chi)$. 
 A {\em non-trivial zero} of $L(s,\chi)$
is a zero of $L(s,\chi)$ such that $0<\Re(s)<1$. (The other zeros 
are called trivial because we know where they are, namely, at negative
integers and, in some cases, also on the line $\Re(s)=0$. In order
to eliminate all zeros on $\Re(s)=0$ outside $s=0$, it suffices to
assume that $\chi$ is {\em primitive}; a primitive character modulo $q$
is one that is not induced by (i.e., not the restriction of) 
any character modulo $d|q$, $d<q$.) 

The Generalized
Riemann Hypothesis for Dirichlet $L$-functions is the statement that, for
every Dirichlet character $\chi$, every non-trivial zero of $L(s,\chi)$
satisfies $\Re(s) = 1/2$. Of course, the Generalized Riemann Hypothesis
(GRH) -- and the Riemann Hypothesis, which is the special case of
$\chi$ trivial --
remains unproven. Thus, if we want to prove unconditional statements, we
need to make do with partial results towards GRH. 
Two kinds of such results have been proven:
\begin{itemize}
\item {\bf Zero-free regions.} Ever since the late nineteenth century
(Hadamard, de la Vall\'ee-Poussin) we have known that there are
hourglass-shaped
regions (more precisely, of the shape $\frac{c}{\log t} \leq \sigma \leq
1 - \frac{c}{\log t}$, where $c$ is a constant and where we write
$s = \sigma + i t$) outside which non-trivial zeros cannot lie. Explicit
values for $c$ are known \cite{MR751161}, \cite{MR2140161}, \cite{Habiba}. 
There is
also the Vinogradov-Korobov region \cite{MR0106205}, \cite{MR0103861},
which is broader asymptotically but narrower in most of
the practical range (see \cite{MR1936814}, however).

\item {\bf Finite verifications of GRH.} It is possible to (ask the
  computer to) prove small, finite fragments of GRH, in the sense of verifying 
that all non-trivial zeros of a given finite set of $L$-functions 
with imaginary part less than some constant $H$ lie on the critical line 
$\Re(s)=1/2$. Such verifications go back to Riemann, who checked the first
few zeros of $\zeta(s)$. Large-scale, rigorous
computer-based verifications are now a possibility.
\end{itemize}

Most work in the literature follows the first alternative, though
\cite{Tao} did use a finite verification of RH (i.e., GRH for the trivial
character). Unfortunately, zero-free regions seem too narrow to be useful
for the ternary Goldbach problem. Thus, we are left with the second
alternative. 

In coordination with the present work, Platt \cite{Plattfresh} verified that
all zeros $s$ of $L$-functions for characters $\chi$ with modulus $q\leq
300000$ satisfying $\Im(s)\leq H_q$ lie on the line $\Re(s)=1/2$, where
\begin{itemize}
\item $H_q = 10^8/q$ for $q$ odd, and
\item $H_q = \max(10^8/q,200+7.5\cdot 10^7/q)$ for $q$ even.
\end{itemize}
This was a medium-large computation, taking a few hundreds of thousands
of core-hours on a parallel computer. It used {\em interval arithmetic} for
the sake of rigor; we will later discuss what this means.

The choice to use a finite verification of GRH, rather than zero-free
regions, had consequences on the manner in which the major and minor
arcs had to be chosen. As we shall see, such a verification can be used
to give very precise bounds on the major arcs, but also forces us to
define them so that they are narrow and their number is constant.
To be precise: the major arcs were defined around rationals $a/q$ with
$q\leq r$, $r=300000$; moreover, as will become clear,
the fact that $H_q$ is finite will
force their width to be bounded by $c_0 r/q x$, where $c_0$ is a constant
(say $c_0=8$). 


\subsection{Estimates of $\widehat{f}(\alpha)$ for $\alpha$ in the major
  arcs}\label{subs:karm}

Recall that we want to estimate sums of the type 
$\widehat{f}(\alpha) = \sum f(n) e(-\alpha n)$, where $f(n)$ is something
like $(\log n) \eta(n/x)$ for $n$ equal to a prime, and $0$ otherwise;
here $\eta:\mathbb{R}\to \mathbb{C}$ is some function of fast decay, 
such as Hardy
and Littlewood's choice, $\eta(t) = e^{-t}$.
 Let us modify this just a little -- we will actually estimate 
\begin{equation}\label{eq:holo}
S_\eta(\alpha,x) = \sum \Lambda(n) e(\alpha n) \eta(n/x),\end{equation}
where $\Lambda$ is the von Mangoldt function (as in (\ref{eq:koj})) .
The use of $\alpha$ rather than $-\alpha$ is just a bow to tradition, as is
the use of the letter $S$ (for ``sum''); however, the use of $\Lambda(n)$
rather than just plain $\log p$ does actually simplify matters.

The function $\eta$ here is sometimes called a {\em smoothing function} or
simply a {\em smoothing}. It will
indeed be helpful for it to be smooth on $(0,\infty)$, but,
in principle, it need not even be continuous. (Vinogradov's work implicitly
uses, in effect, the ``brutal truncation'' $1_{\lbrack 0,1\rbrack}(t)$,
defined to be $1$ when $t\in \lbrack 0,1\rbrack$ and $0$ otherwise; that
would be fine for the minor arcs, but, as it will become clear, it is a bad idea as
far as the major arcs are concerned.)

Assume $\alpha$ is on a major arc, meaning that we can write 
$\alpha = a/q + \delta/x$ for some $a/q$ ($q$ small) and some $\delta$ 
(with $|\delta|$ small). We can write
 $S_\eta(\alpha,x)$ as a linear combination 
\begin{equation}\label{eq:sorrow}
S_\eta(\alpha,x) = \sum_\chi c_\chi
S_{\eta,\chi}\left(\frac{\delta}{x},x\right)
+ \text{tiny error term},\end{equation} where 
\begin{equation}\label{eq:battle}
S_{\eta,\chi}\left(\frac{\delta}{x},x\right) = \sum \Lambda(n) \chi(n)
e(\delta n/x) \eta(n/x).\end{equation}
In (\ref{eq:sorrow}),
 $\chi$ runs over primitive Dirichlet characters of moduli $d|q$, and
$c_\chi$ is small ($|c_\chi|\leq \sqrt{d}/\phi(q)$).

Why are we expressing the sums $S_\eta(\alpha,x)$ 
in terms of the sums $S_{\eta,\chi}(\delta/x,x)$,
which look more complicated?
The argument has become $\delta/x$, whereas before it was $\alpha$. 
Here $\delta$ is relatively small 
-- smaller than the constant $c_0 r$, in our setup. 
In other words, $e(\delta n/x)$ will go around the circle a bounded number of 
times as $n$ goes from $1$ up to a constant times $x$ 
(by which time $\eta(n/x)$ has become small, because $\eta$ is of fast decay). 
This makes the sums much easier to estimate.

To estimate the sums $S_{\eta,\chi}$, we will use $L$-functions, together
with one of the most common tools of analytic number theory, the Mellin
transform.
This transform is essentially a Laplace transform with a change of
variables, and a Laplace transform, in turn, is a Fourier transform taken
on a vertical line in the complex plane. For $f$ of fast enough decay,
the Mellin transform $F=Mf$ of $f$ is given by
\[F(s) = \int_0^\infty f(t) t^s \frac{dt}{t};\]
we can express $f$ in terms of $F$ by the {\em Mellin inversion formula}
\[f(t) = \frac{1}{2\pi i} \int_{\sigma-i\infty}^{\sigma+i\infty} F(s) t^{-s} ds
\]
for any $\sigma$ within an interval. We can thus express
$e(\delta t) \eta(t)$ in terms of its Mellin transform $F_\delta$ and
then use (\ref{eq:koj}) to express $S_{\eta,\chi}$ in terms of $F_\delta$
and $L'(s,\chi)/L(s,\chi)$; shifting the integral in the
Mellin inversion formula to the left, we obtain what is known in analytic
number theory as an {\em explicit formula}: 
\[S_{\eta,\chi}(\delta/x,x) = \left\lbrack \widehat{\eta}(
-\delta) x\right\rbrack 
 - \sum_\rho F_\delta(\rho) x^\rho + \text{tiny error term}.\]
Here the term between brackets appears only for $\chi$ trivial. In the sum, 
$\rho$ goes over all non-trivial zeros of $L(s,\chi)$, and $F_\delta$ 
is the Mellin transform of $e(\delta t) \eta(t)$. (The tiny error term
comes from a sum over the trivial zeros of $L(s,\chi)$.)
We will obtain the estimate we desire
 if we manage to show that the sum over $\rho$ is small.

The point is this: if we verify GRH for $L(s,\chi)$ up to imaginary part
$H$, i.e., if we check that
all zeroes $\rho$ of $L(s,\chi)$ 
with $|\Im(\rho)|\leq H$ satisfy $\Re(\rho)=1/2$, we have
$|x^\rho| = \sqrt{x}$. In other words, $x^\rho$ 
is very small (compared to $x$). However, for any $\rho$ whose imaginary part has absolute value greater than $H$, we know next to nothing about its real part, other than $0\leq \Re(\rho)\leq 1$. (Zero-free regions
are notoriously weak for $\Im(\rho)$ large; we will not use them.) 
Hence, our only chance is to make sure that $F_\delta(\rho)$ is very small when $|\Im(\rho)|\geq H$. 

This has to be true for both $\delta$ very small 
(including the case $\delta=0$) and for $\delta$ not so small ($|\delta|$
up to $c_0 r/q$, which can be large because 
 $r$ is a large constant). How can we choose $\eta$ so that
$F_\delta(\rho)$ is very small in both cases for $\tau=\Im(\rho)$ large?

The method of {\em stationary phase} is useful as an exploratory tool here.
In brief, it suggests (and can sometimes prove) that the main 
contribution to the integral 
\begin{equation}\label{eq:narimsi}
F_\delta(t) = \int_0^\infty e(\delta t) \eta(t) t^s \frac{dt}{t}\end{equation}
can be found where 
the phase of the integrand has derivative $0$. 
This happens when $t= -\tau/2\pi \delta$ (for $\sgn(\tau)\ne \sgn(\delta)$);
the contribution is then a moderate factor times $\eta(-\tau/2\pi \delta)$.
In other words, if $\sgn(\tau)\ne \sgn(\delta)$ and $\delta$ is not too small
($|\delta|\geq 8$, say),
$F_\delta(\sigma + i\tau)$ behaves like $\eta(-\tau/2\pi \delta)$;
if $\delta$ is small ($|\delta|<8$), then $F_\delta$ behaves like $F_0$,
which is the Mellin transform $M\eta$ of $\eta$. Here is our goal, then: the decay of
$\eta(t)$ as $|t|\to \infty$ should be as fast as possible, and the
decay of the transform $M\eta(\sigma+ i \tau)$ should also be as fast as possible.

This is a classical dilemma, often called the uncertainty principle because
it is the mathematical fact underlying the physical principle of the same name:
you cannot have a function  
$\eta$ that decreases extremely rapidly and whose Fourier
transform (or, in this case its Mellin transform) also decays extremely 
rapidly.

What does ``extremely rapidly'' mean here? It means (as Hardy himself proved)
``faster than any exponential $e^{- C t}$''. Thus, Hardy and Littlewood's
choice $\eta(t) = e^{-t}$ seems essentially optimal at first sight.

However, it is not optimal. We can choose $\eta$ so that $M\eta$
decreases exponentially (with a constant $C$ somewhat worse than for
$\eta(t)=e^{-t}$),
but $\eta$ decreases faster than exponentially. This is a particularly
appealing possibility because it is $t/|\delta|$, and not so much $t$,
that risks being fairly small. (To be explicit: say we check GRH for
characters of modulus $q$ up to $H_q\sim 50 \cdot c_0 r/q \geq 50 |\delta|$.
Then we only know that $|\tau/2\pi\delta| \gtrsim 8$. So, for
$\eta(t) = e^{-t}$, $\eta(-\tau/2\pi \delta)$ may be as large as $e^{-8}$,
which is not negligible. Indeed, since this term will be multiplied
later by other terms, $e^{-8}$ is simply not small enough. 
On the other hand,
we can assume that $H_q\geq 200$ (say), and so $M\eta(s) \sim e^{-(\pi/2)
|\tau|}$ is completely negligible, and will remain negligible even
if we replace $\pi/2$ by a somewhat smaller constant.)

We shall take $\eta(t) = e^{-t^2/2}$ (that is, the Gaussian). 
This is not the only 
possible choice, but it is in some sense natural. It is easy to show that 
the Mellin transform $F_\delta$ for $\eta(t) = e^{-t^2/2}$ is a multiple of what
is called a {\em parabolic cylinder function} $U(a,z)$ 
with imaginary values for $z$. There are plenty of 
estimates on parabolic cylinder functions in the literature -- but mostly
for $a$ and $z$ real, in part because that is one
 of the cases occuring most often in applications. There
are some asymptotic expansions and estimates for $U(a,z)$, $a$, $z$,
general, due to Olver \cite{MR0094496}, 
\cite{MR0109898}, \cite{MR0131580}, \cite{MR0185350},
but unfortunately they come without fully explicit error terms for $a$
and $z$ within our range of interest. (The same holds for \cite{MR1993339}.)

In the end, I derived bounds for $F_\delta$ using the {\em saddle-point method}.
(The method of stationary phase, which we used to choose $\eta$, seems to lead
to error terms that are too large.) The saddle-point method
consists, in brief, in changing the contour of an integral to be bounded
(in this case, (\ref{eq:narimsi})) so as to minimize the maximum of
the integrand, and so as to go as quickly as possible through the point
at which the maximum is reached. (To use a metaphor in \cite{MR671583}:
find the lowest mountain pass and descend from it as quickly as possible.)
The interesting part here (as, it seems, in other applications of the method)
is to find a contour satisfying these conditions while leading to an integral
that can be estimated relatively cleanly. (The use of rigorous numerics --
to give bounds on extrema and series expansions,
rather than to perform integration -- was also helpful here.)

For $s = \sigma + i\tau$ with $\sigma\in \lbrack 0,1\rbrack$ and
$|\tau|\geq \max(100,4\pi^2 |\delta|)$, we obtain that the Mellin transform
$F_\delta$ of $\eta(t) e(\delta t)$ with $\eta(t) = e^{-t^2/2}$ satisfies
\begin{equation}\label{eq:sturmo}
|F_\delta(s)| + |F_\delta(1-s)| \leq 4.226 \cdot \begin{cases}
e^{-0.1065 \left(\frac{\tau}{\pi \delta}\right)^2} &\text{if $|\tau|< \frac{3}{2}
(\pi \delta)^2$,}\\ e^{-0.1598 |\tau|} &\text{if $|\tau|\geq \frac{3}{2}
(\pi \delta)^2$.}\end{cases}\end{equation}
Similar bounds hold for $\sigma$ in other ranges, thus giving us (similar)
estimates for the Mellin transform $F_\delta$ for $\eta(t) = t^k e^{-t^2/2}$
and $\sigma$ in the critical range $\lbrack 0,1\rbrack$. 

A moment's thought shows that we can also use (\ref{eq:sturmo}) to deal
with the Mellin transform of $\eta(t) e(\delta t)$ for any function of the
form $\eta(t) = e^{-t^2/2} g(t)$ (or, more generally, $\eta(t) = t^k e^{-t^2/2}
g(t)$), where $g(t)$ is any {\em band-limited function}. By a band-limited
function, we could mean a function whose Fourier transform is compactly
supported; while that is a plausible choice, it turns out to be better to
work with functions that are band-limited with respect to the Mellin transform
-- in the sense of being of the form
\[g(t) = \int_{-R}^R h(r) t^{-ir} dr,\]
where $h:\mathbb{R}\to \mathbb{C}$ is supported on a compact interval 
$\lbrack -R,R\rbrack$, with $R$ not too large (say $R=200$). What happens
is that the Mellin transform of the product $e^{-t^2/2} g(t) e(\delta t)$ 
is a convolution of the Mellin transform $F_\delta(s)$ 
of $e^{-t^2/2} e(\delta t)$ 
(estimated in (\ref{eq:sturmo})) and that of $g(t)$ (supported in 
$\lbrack -R,R\rbrack$); the effect of the convolution is just to delay
decay of $F_\delta(s)$ by, at most, a shift by $y\mapsto y-R$.

There remains to do one thing, namely, to derive an explicit formula
general enough to work with all the weights $\eta(t)$ we have
discussed and some we will discuss later, while being
 also completely explicit, and free of any integrals that 
may be tedious to evaluate. Once that is done, and once we consider the
input provided by Platt's finite verification of GRH up to $H_q$, we
obtain simple bounds for different weights.

For $\eta(t)= e^{-t^2/2}$, $x\geq 10^8$, $\chi$ a primitive character
of modulus $q\leq r = 300000$, and any $\delta\in \mathbb{R}$ with
$|\delta|\leq 4r/q$, we obtain
\begin{equation}\label{eq:frank}
S_{\eta,\chi}\left(\frac{\delta}{x}, x\right)
= I_{q=1} \cdot \widehat{\eta}(-\delta) x
+ E\cdot x,\end{equation}
where $I_{q=1}=1$ if $q=1$, $I_{q=1}=0$ if $q\ne 1$,
and
\begin{equation}\label{eq:lynno}
|E|\leq
5.281 \cdot 10^{-22} + \frac{1}{\sqrt{x}}
\left( \frac{650400}{\sqrt{q}} + 112\right).\end{equation}
Here $\widehat{\eta}$ stands for the Fourier transform from $\mathbb{R}$
to $\mathbb{R}$ normalized as follows: $\widehat{\eta}(t) = 
\int_{-\infty}^\infty e(-xt) \eta(x) dx$. Thus, $\widehat{\eta}(-\delta)$
is just $\sqrt{2\pi} e^{-2 \pi^2 \delta^2}$ (self-duality of the
Gaussian).

This is one of the main results of \cite{HelfMaj}. Similar bounds are
also proven there for $\eta(t) = t^2 e^{-t^2/2}$, as well as for a weight
of type $\eta(t) = t e^{-t^2/2} g(t)$, where $g(t)$ is a band-limited
function, and also for a weight $\eta$ defined by a multiplicative convolution. 
The conditions on $q$ ($q\leq r = 300000$) and $\delta$ are what we
expected from the outset.

Thus concludes our treatment of the major arcs. This is arguably the easiest
part of the proof; it was actually what I left for the end, as I was fairly
confident it would work out. Minor-arc estimates are more delicate; 
let us now examine them.

\section{The minor arcs $\mathfrak{m}$}\label{sec:melusa}

\subsection{Qualitative goals and main ideas}
What kind of bounds do we need? What is there in the literature?

We wish to obtain upper bounds on $|S_\eta(\alpha,x)|$ for some weight $\eta$
and any $\alpha\in \mathbb{R}/\mathbb{Z}$ not very close to a rational
with small denominator. Every $\alpha$ is close to some rational $a/q$;
what we are looking for is a bound on $|S_\eta(\alpha,x)|$ that
decreases rapidly when $q$ increases.

Moreover, we want our bound to decrease rapidly when $\delta$ increases,
where $\alpha = a/q+\delta/x$. In fact, the main terms in our bound will
be decreasing functions of $\max(1,|\delta|/8)\cdot q$.
(Let us write $\delta_0 = \max(2,|\delta|/4)$ from now on.)
This will allow our bound to be good enough outside 
 narrow major arcs, which will get narrower and narrower as $q$ increases
-- that is, precisely the kind of major arcs we were presupposing in 
our major-arc bounds.

It would be possible to work with narrow major arcs that
become narrower as $q$ increases simply by allowing $q$ to be very large 
(close to $x$), and assigning each angle to the fraction closest to it.
This is, in fact, the common procedure. 
However, this makes matters more difficult, in that we would have to 
minimize at the same time the factors in front of terms $x/q$, $x/\sqrt{q}$,
etc., and those in front of terms $q$, $\sqrt{q x}$, and so on. 
(These terms are being compared to the trivial bound $x$.)
Instead, 
we choose to strive for a direct dependence on $\delta$ throughout; this
will allow us to cap $q$ at a much lower level, thus making terms such as
$q$ and $\sqrt{q x}$ negligible. (This choice has been taken
elsewhere in applications of the circle method, but, strangely, seems absent
from previous work on the ternary Goldbach conjecture.)

How good must our bounds be? Since the major-arc bounds are valid only
for $q\leq r=300000$ and $|\delta|\leq 4r/q$, we cannot afford even
a single factor of $\log x$ (or any other function tending to $\infty$
as $x\to \infty$) in front of terms such as $x/\sqrt{q |\delta_0|}$:
a factor like that would make the term larger than the trivial bound $x$
for $q |\delta_0|$ equal to a constant ($r$, say) and $x$ very large. 
Apparently, there was no such ``log-free bound''
with explicit constants in the literature, even though
such bounds
 were considered to be in principle feasible, and even though previous
work (\cite{MR813837}, \cite{MR1399341}, \cite{MR1803131}, \cite{Tao})
had gradually decreased the number of factors of $\log x$.
(In limited ranges for $q$, there were log-free bounds without
explicit constants; see \cite{MR1399341}, \cite{MR2607306}. The estimate in
\cite[Thm. 2a, 2b]{MR0062183} was almost log-free, but not quite.
There were also bounds \cite{MR1215269}, \cite{MR2776653}
that used $L$-functions, and thus were not
really useful in a truly minor-arc regime.)

It also seemed clear that a main bound proportional to $(\log q)^2 x/\sqrt{q}$
(as in \cite{Tao}) was too large. At the same time, it was not really
necessary to reach a bound of the best possible form that could
be found through Vinogradov's basic approach, namely
\begin{equation}\label{eq:astora}
|S_\eta(\alpha,x)|\leq C \frac{x \sqrt{q}}{\phi(q)}.\end{equation}
Such a bound had been proven by Ramar\'e \cite{MR2607306} for $q$
in a limited range and $C$ non-explicit; later, in \cite{Ramlater}
-- which postdates the first version of \cite{Helf} --
Ramar\'e  broadened the
range to $q\leq x^{1/48}$ and gave an explicit value for $C$, namely,
$C=13000$. Such a bound is a notable achievement, but, unfortunately, it
is not useful for our purposes. Rather, we will aim at a bound whose
main term is bounded by a constant around $1$ times 
$x (\log \delta_0 q)/\sqrt{\delta_0 \phi(q)}$; this is slightly worse
asymptotically than (\ref{eq:astora}), but it is much better in the delicate
range of $\delta_0 q \sim 300000$, and in fact for a much wider range as well.

\begin{center}
* * *
\end{center}

We see that we have several tasks. One of them is the removal of
logarithms: we cannot afford a single factor of $\log x$, and, in practice,
we can afford at most one factor of $\log q$. Removing logarithms will
be possible in part because of the use of efficient techniques (the large
sieve for sequences with prime support) but also because we will be able to find cancellation
at several places in sums coming from a combinatorial identity (namely, 
Vaughan's identity). The task of finding cancellation is particularly delicate
because we cannot afford large constants or, for that matter,
 statements valid only for large
$x$. (Bounding a sum such as $\sum_n \mu(n)$ efficiently (where $\mu$
is the M\"obius function)
is harder than estimating a sum such as $\sum_n \Lambda(n)$ equally efficiently, even though we are used to thinking of the two problems as equivalent.)

We have said that our bounds will improve as $|\delta|$ increases. This
dependence on $\delta$ will be secured in different ways at different
places. Sometimes $\delta$ will appear as an argument, as in
$\widehat{\eta}(-\delta)$; for $\eta$ piecewise continuous with $\eta' \in
L_1$, we know that
$|\widehat{\eta}(t)| \to 0$ as $|t|\to \infty$. Sometimes we will
obtain a dependence on $\delta$ by using several different rational 
approximations to the same $\alpha \in \mathbb{R}$. Lastly, 
 we will obtain a good dependence on $\delta$ in bilinear
sums by supplying a scattered input to a large sieve.

If there is a main moral to the argument, it lies in the 
close relation between the circle method and
the large sieve. The circle method rests on the
estimation of an integral involving a Fourier transform $\widehat{f}:
\mathbb{R}/\mathbb{Z}\to \mathbb{C}$; as we will later see, this leads
naturally to estimating the $\ell_2$-norm of $\widehat{f}$ on subsets 
(namely, unions of arcs) of the circle $\mathbb{R}/\mathbb{Z}$. The large
sieve can be seen as an approximate discrete version of Plancherel's identity,
which states that $|\widehat{f}|_2 = |f|_2$. 

Both in this section and in \S \ref{sec:putall}, 
we shall use the large sieve in part so as to use the fact that some of the 
functions we work with have prime support, i.e., are non-zero only
on prime numbers. There are
ways to use prime support to improve the output of the large sieve.
 In \S \ref{sec:putall}, these techniques will be refined and
then translated to the context of the circle method,
where $f$ has (essentially) prime support and $|\widehat{f}|^2$ must be 
integrated over unions of arcs. (This allows us to remove a logarithm.)
The main point is that the large sieve is not being used as a black
box; rather, we can adapt ideas from (say) the large-sieve context and apply
them to the circle method.

Lastly, there are the benefits of a continuous $\eta$. Hardy and
Littlewood already used a continuous $\eta$; this was abandoned by Vinogradov, 
presumably for the sake of simplicity. The idea that smooth weights $\eta$
can be superior to sharp truncations is now commonplace. As we shall see,
using a continuous $\eta$ is helpful in the minor-arcs regime, but not
as crucial there as for the major arcs. We will not use a smooth $\eta$; we 
will prove our estimates for any continuous $\eta$ that is piecewise $C_1$,
and then, towards the end, we will choose to 
use the same weight $\eta=\eta_2$ as in \cite{Tao}, in part because it
has compact support, and in part for the sake of comparison. The moral here
is not quite the common dictum ``always smooth'', but rather that different kinds of smoothing can 
be appropriate for different tasks; in the end, we will show how
to coordinate different smoothing functions $\eta$. 

There are other ideas involved; for instance, some of Vinogradov's lemmas  
are improved. Let us now go into some of the details.
 
\subsection{Combinatorial identities}

Generally, since Vinogradov, a 
treatment of the minor arcs starts with a combinatorial identity expressing
$\Lambda(n)$ (or the characteristic function of the primes) as a sum of
two or more convolutions. (In this section, by a convolution $f\ast g$, we will 
mean the {\em Dirichlet convolution} $(f\ast g)(n) = \sum_{d|n} f(d) g(n/d)$,
i.e., the multiplicative convolution on the semigroup of positive integers.)

In some sense, the archetypical identity is \[\Lambda = \mu\ast \log,\]
but it will not usually do: the contribution of $\mu(d) \log(n/d)$ with
$d$ close to $n$ is too difficult
to estimate precisely. There are alternatives: for example, there is
Selberg's identity
\begin{equation}\Lambda(n) \log n = \mu\ast \log^2 - \Lambda\ast \Lambda,
\end{equation}
or the generalization of this to $\Lambda(n) (\log n)^k = \mu\ast
\log^{k+1} - \dotsc$ (Bomberi-Selberg), 
used in Bomberi's strengthening of the
Erd\H{o}s-Selberg proof of the prime number theorem.
Another useful (and very simple) identity was that used by
Daboussi's \cite{MR1803131}; see also
\cite{MR1399341}, which gives explicit estimates of sums over primes.

The proof of Vinogradov's three-prime result was simplified substantially 
in \cite{MR0498434} by the introduction of {\em Vaughan's identity}:
\begin{equation}\label{eq:vaughan}
\Lambda(n) = \mu_{\leq U} \ast \log - \Lambda_{\leq V} \ast \mu_{\leq U} \ast 1
+ 1 \ast \mu_{>U} \ast \Lambda_{>V} + \Lambda_{\leq V},\end{equation}
where we are using the notation
\[f_{\leq W} = \begin{cases} f(n) &\text{if $n\leq W$,}\\ 0 
&\text{if $n>W$,}\end{cases}\;\;\;\;\;\;\;
f_{>W} = \begin{cases} 0 &\text{if $n\leq W$,}\\ f(n) 
&\text{if $n>W$.}\end{cases}\]
Of the resulting sums ($\sum_n (\mu_{\leq U}\ast \log)(n) e(\alpha n) \eta(n/x)$,
etc.), the first three
 are said to be of {\em type I}, {\em type I} (again) and {\em type II};
the last sum, $\sum_{n\leq V} \Lambda(n)$, is negligible.

One of the advantages of Vaughan's identity is its flexibility: we can
set $U$ and $V$ to whatever values we wish. Its main disadvantage is that
it is not ``log-free'', in that it seems to impose the loss of two factors of
$\log x$: if we sum each side of (\ref{eq:vaughan}) from $1$ to $x$,
we obtain $\sum_{n\leq x} \Lambda(n) \sim x$ on the left side, whereas,
if we bound the sum on the right side without the use of cancellation,
we obtain a bound of $x (\log x)^2$. Of course, we will obtain some
 cancellation from the phase $e(\alpha n)$; still, even if this gives
us a factor of, say, $1/\sqrt{q}$, we will get a bound of 
$x (\log x)^2/\sqrt{q}$,
which is worse than the trivial bound $x$ for $q$ bounded and $x$ large.
Since we want a bound that is useful for all $q$ larger than the constant $r$
and all $x$ larger than a constant, this will not do. 

As was pointed out in \cite{Tao}, it is possible to get a factor of
$(\log q)^2$ instead of a factor of $(\log x)^2$ in the type II sums
by setting $U$ and $V$ appropriately. Unfortunately, a factor of $(\log q)^2$
is still too large in practice, and there is also the issue of factors of
$\log x$ in type I sums.

Vinogradov had already managed to get an essentially log-free result
(by a rather difficult procedure) in \cite[Ch. IX]{MR0062183}. 
The result in \cite{MR1399341} is log-free. Unfortunately, the explicit result
in \cite{MR1803131} -- the study of which encouraged me at the beginning of 
the project -- is not. For a while, I worked with the Bombieri-Selberg
identity with $k=2$. Ramar\'e obtained a log-free bound in \cite{MR2607306}
using the Diamond-Steinig identity, which is related to Bombieri-Selberg.

In the end, I decided to use Vaughan's identity.
This posed a challenge: to obtain cancellation in Vaughan's identity 
at every possible step, beyond the cancellation given by the phase 
$e(\alpha n)$. (The presence of a phase, in fact, makes the task of
getting cancellation from the identity more complicated.)
The removal of logarithms will be one of
our main tasks in what follows. It is clear that the presence of the 
M\"{o}bius function $\mu$ should give, in principle, some cancellation; we
will show how to use it to obtain as much cancellation as we need -- 
with good constants, and not just asymptotically.

\subsection{Type I sums}
There are two type I sums, namely,
\begin{equation}\label{eq:rada1}
\sum_{m\leq U} \mu(m) \sum_n (\log n) e(\alpha m n) \eta\left(\frac{m n}{x}\right)
\end{equation}
and
\begin{equation}\label{eq:rada2}
\sum_{v\leq V} \Lambda(v) \sum_{u\leq U} \mu(u) \sum_n
e(\alpha v u n) \eta\left(\frac{v u n}{x}\right).\end{equation}
In either case, $\alpha = a/q + \delta/x$, where $q$ is larger than a constant
$r$ and $|\delta/x|\leq 1/q Q_0$ for some $Q_0>\max(q,\sqrt{x})$.
For the purposes of this exposition, we will set it as our task to 
estimate the slightly simpler
sum
\begin{equation}\label{eq:sadun}
\sum_{m\leq D} \mu(m) \sum_n e(\alpha m n) \eta\left(\frac{m n}{x}\right),\end{equation}
where $D$ can be $U$ or $UV$ or something else less than $x$.

Why can we consider this simpler sum without omitting anything essential?
It is clear that (\ref{eq:rada1}) is of the same kind as (\ref{eq:sadun}).
The inner double sum in (\ref{eq:rada2}) is just (\ref{eq:sadun})
with $\alpha v$ instead of $\alpha$; this enables us to estimate
(\ref{eq:rada2}) by means of (\ref{eq:sadun}) for $q$ small, i.e., the
more delicate case. If $q$ is not small, then the approximation
$\alpha v \sim a v/q$ may not be accurate enough. In that case, we
collapse the two outer sums in (\ref{eq:rada2}) into 
a sum $\sum_n (\Lambda_{\leq V}\ast \mu_{\leq U})(n)$, and treat all of
(\ref{eq:rada2}) much as we will treat (\ref{eq:sadun}); since $q$ is not
small, we can afford to bound $(\Lambda_{\leq V}\ast \mu_{\leq U})(n)$ trivially
(by $\log n$) in the less sensitive terms.

Let us first outline Vinogradov's procedure for bounding type I sums.
Just by summing a geometric series, we get
 \begin{equation}\label{eq:salome}
\left|\sum_{n\leq N} e(\alpha n) \right| \leq
\min\left(N,\frac{c}{\{\alpha\}}\right),\end{equation}
where $c$ is a constant and
 $\{\alpha\}$ is the distance from $\alpha$ to the nearest integer.
Vinogradov splits the outer sum in (\ref{eq:sadun}) into sums of length
$q$. When $m$ runs on an interval of length $q$, the angle $a m/q$ runs
through all fractions of the form $b/q$; due to the error $\delta/x$,
$\alpha m$ could be close to $0$ for two values of $n$, but otherwise
$\{\alpha m\}$ takes values bounded below by $1/q$ (twice), $2/q$ (twice),
$3/q$ (twice), etc. Thus
\begin{equation}\label{eq:kolt}\left|\sum_{y<m\leq y+q} \mu(m) \sum_{n\leq N} e(\alpha m n)\right|
\leq \sum_{y<m\leq y+q} \left|\sum_{n\leq N} e(\alpha m n)\right| \leq
\frac{2 N}{m} + 2 c q \log e q\end{equation}
for any $y\geq 0$.

There are several ways to improve this. One is simply to estimate
the inner sum more precisely; this was already done in \cite{MR1803131}.
One can also define a smoothing function $\eta$, as in (\ref{eq:sadun});
it is easy to get
\[\left|\sum_{n\leq N} e(\alpha n) \eta\left(\frac{n}{x}\right)\right|\leq
\min\left(x |\eta|_1 + \frac{|\eta'|_1}{2}, \frac{|\eta'|_1}{2
|\sin(\pi \alpha)|}, \frac{|\widehat{\eta''}|_\infty}{4 
x (\sin \pi \alpha)^2}\right).\]
Except for the third term, this is as in \cite{Tao}. We could also choose
carefully which bound to use for each $m$; surprisingly, this gives
an improvement -- in fact, an important one, for $m$ large. 
However, even with these improvements, we still
have a term proportional to $N/m$ as in (\ref{eq:kolt}), and this contributes
about $(x \log x)/q$ to the sum (\ref{eq:sadun}), thus giving us an estimate
that is not log-free.

What we have to do, naturally, is to take out the terms with $q|m$ for
$m$ small. (If $m$ is large, then those may not be the terms for which
$m\alpha$ is close to $0$; we will later see what to do.) For $y+q\leq Q/2$,
$|\alpha - a/q|\leq 1/q Q$, 
we get that
\begin{equation}\label{eq:cami}\mathop{\sum_{y<m\leq y+q}}_{q\nmid m} 
 \min\left(A, \frac{B}{|\sin \pi \alpha n|},
\frac{C}{|\sin \pi \alpha n|^2}\right)\end{equation}
is at most
\begin{equation}\label{eq:crims}
\min\left(\frac{20}{3 \pi^2} C q^2, 2 A + \frac{4 q}{\pi} \sqrt{AC},
\frac{2 B q}{\pi} \max\left(2, \log \frac{C e^3 q}{B \pi}\right)\right).
\end{equation}
This is satisfactory. We are left with all the terms $m\leq M = \min(D,Q/2)$
with $q|m$ -- and also with all the terms $Q/2 < m\leq D$.
For $m\leq M$ divisible by $q$, we can estimate (as opposed
to just bound from above) 
the inner sum in (\ref{eq:sadun}) by the Poisson summation
formula, and then sum over $m$, but without taking absolute values; writing
$m=aq$, we get a main term 
\begin{equation}\label{eq:jut}
\frac{x \mu(q)}{q} \cdot \widehat{\eta}(-\delta) \cdot
\mathop{\sum_{a\leq M/q}}_{(a,q)=1} \frac{\mu(a)}{a},\end{equation}
where $(a,q)$ stands for the greatest common divisor of $a$ and $q$.

It is clear that we have to get cancellation over $\mu$ here.
There is an elegant elementary argument \cite{MR1401709} showing that
the absolute value of the sum in (\ref{eq:jut}) is at most $1$.
We need to gain one more log, however. Ramar\'e \cite{Ramsev} helpfully
furnished the following bound:
\begin{equation}\label{eq:hujt}
\left|\mathop{\sum_{a\leq x}}_{(a,q)=1} \frac{\mu(a)}{a}\right|\leq
\frac{4}{5} \frac{q}{\phi(q)} \frac{1}{\log x/q}\end{equation}
for $q\leq x$. (Cf. \cite{MR1378588}, \cite{ElMarraki})
This is neither trivial nor elementary.\footnote{The current state of knowledge may seem surprising: after all, we expect
nearly square-root cancellation ($|\sum_{n\leq x} \mu(n)/n|\leq \sqrt{2/x}$
holds for all real $0<x\leq 10^{12}$; see also the stronger bound 
\cite{MR1259423}).
The classical zero-free region of the
Riemann zeta function ought to give a factor of $e^{-\sqrt{(\log x)/c}}$,
which looks much better than $1/\log x$. What happens is that (a) such a factor
is not actually much better than $1/\log x$ for $x\sim 10^{30}$, say; (b)
estimating sums involving the M\"obius function by means of an explicit formula
is harder than estimating sums involving $\Lambda(n)$:
the residues of $1/\zeta(s)$ at the non-trivial zeros of $s$ come into play.
As a result, getting non-trivial
 explicit results on sums of $\mu(n)$ is more difficult than
one would naively expect from the quality of classical effective results. See 
\cite{RamEtatLieux} for a survey of explicit bounds.} We are, so to speak, allowed
to use non-elementary means (that is, methods based on $L$-functions)
because the only $L$-function we need to use here is the Riemann zeta function.

What shall we do for $m>Q/2$? We can always give a bound
\begin{equation}\label{eq:tuffblo}
\sum_{y<m\leq y+q}
 \min\left(A, 
\frac{C}{|\sin \pi \alpha n|^2}\right) \leq 3 A + \frac{4 q}{\pi} \sqrt{A C} 
\end{equation}
for $y$ arbitrary; since $A C$ will be of constant size, $(4 q/\pi) \sqrt{A C}$
is pleasant enough, but the contribution of $3 A \sim 3 |\eta|_1 x/y$
is nasty (it adds a multiple of $(x \log x)/q$ to the total) and seems
 unavoidable: the values of $m$ for which
$\alpha m$ is close to $0$ no longer correspond to the congruence class
$m\equiv 0 \mo q$, and thus cannot be taken out. 

The solution is to switch approximations. (The idea of using different
approximations to the same $\alpha$ is neither new nor recent in the general
context of the circle method: see \cite[\S 2.8, Ex. 2]{MR1435742}. What
may be new is its use to clear a hurdle in type I sums.) What does this mean?
If $\alpha$ 
were exactly, or almost exactly, $a/q$, then there would be no other very
good approximations in a reasonable range. However, note that we can
{\em define} $Q = \lfloor x/ |\delta q|\rfloor$ for 
$\alpha = a/q + \delta/x$, and still have
$|\alpha - a/q|\leq 1/q Q$. If $\delta$ is very small, $Q$ will be
larger than $2 D$, and there will be no terms with $Q/2 < m \leq D$ to
worry about.

What happens if $\delta$ is not very small? We know that, for any $Q'$, there is
an approximation $a'/q'$ to $\alpha$ with $|\alpha - a'/q'|\leq 1/q' Q'$
and $q'\leq Q'$. However, for $Q' > Q$, we know that $a'/q'$ cannot equal
$a/q$: by the definition of $Q$, the approximation $a/q$ is not good
enough, i.e., $|\alpha - a/q|\leq 1/q Q'$ does not hold.
Since $a/q\ne a'/q'$, we see that $|a/q - a'/q'|\geq 1/q q'$, and this implies that $q'\geq
(\epsilon/(1+\epsilon)) Q$. 

Thus, for $m>Q/2$, the solution is to apply (\ref{eq:tuffblo}) with
$a'/q'$ instead of $a/q$. The contribution of $A$ fades into
insignificance: for the first sum over a range $y<m\leq y+q'$, $y\geq Q/2$,
it contributes at most $x/(Q/2)$, and all the other contributions of $A$
sum up to at most a constant times $(x \log x)/q'$.

Proceeding in this way, we obtain a total bound for
(\ref{eq:sadun}) whose main terms are proportional to
\begin{equation}\label{eq:marxophone}\frac{1}{\phi(q)} \frac{x}{\log \frac{x}{q} }
\min\left(1, \frac{1}{\delta^2}\right),
\;\;\; \frac{2}{\pi} \sqrt{|\widehat{\eta''}|_\infty} \cdot D\;\;\; 
\text{and} \;\;\;\;
q \log \max\left(\frac{D}{q},q\right),
\end{equation}
with good, explicit constants. The first term -- usually the largest one --
is precisely what we needed: it is proportional to $(1/\phi(q)) x/\log x$
for $q$ small, and decreases rapidly as $|\delta|$ increases.

\subsection{Type II, or bilinear, sums}
We must now bound
\[S = \sum_m (1\ast \mu_{>U})(m) \sum_{n>V} \Lambda(n) e(\alpha m n) \eta(m n/x).
\]
At this point it is convenient to assume that $\eta$ is the Mellin
convolution of two functions. The
{\em multiplicative} or {\em Mellin convolution} on
$\mathbb{R}^+$ is defined by
\[(\eta_0\ast_M \eta_1)(t) = \int_0^\infty \eta_0(r)
\eta_1\left(\frac{t}{r}\right) \frac{dr}{r}.\]
Tao \cite{Tao} takes $\eta = \eta_2 = \eta_1 \ast_M \eta_1$, where 
$\eta_1$ is a brutal truncation, viz.,
the function taking the value $2$ on $\lbrack 1/2,1\rbrack$
and $0$ elsewhere. We take the same $\eta_2$, in part for comparison
purposes, and in part because this will allow us to use off-the-shelf 
estimates on the large sieve. (Brutal truncations are rarely optimal in
principle,
but, as they are very common, results for them have been carefully optimized
in the literature.) Clearly
\begin{equation}\label{eq:salom}
S = \int_V^{X/U} \sum_m \left(\mathop{\sum_{d>U}}_{d|m}
\mu(d)\right) \eta_1\left(\frac{m}{x/W}\right) \cdot \sum_{n\geq V}
\Lambda(n) e(\alpha m n) \eta_1\left(\frac{n}{W}\right) \frac{dW}{W}. \end{equation}
By Cauchy-Schwarz, the integrand is at most $\sqrt{S_1(U,W) S_2(V,W)}$, where
\begin{equation}\label{eq:katar}\begin{aligned}
S_1(U,W) &= \sum_{\frac{x}{2 W} < m\leq \frac{x}{W}} \left|\mathop{\sum_{d>U}}_{d|m} \mu(d)\right|^2,
\\ S_2(V,W) &= \sum_{\frac{x}{2 W} \leq m\leq \frac{x}{W}} \left|
\sum_{\max\left(V,\frac{W}{2}\right) \leq n\leq W}
\Lambda(n) e(\alpha m n)\right|^2.\end{aligned}\end{equation}

We must bound $S_1(U,W)$ by a constant times $x/W$. We are able to do this
-- with a good constant. (A careless bound would have given a multiple of
$(x/U) \log^3 (x/U)$, which is much too large.) First,
we reduce $S_1(W)$ to an expression involving an integral of
\begin{equation}\label{eq:rotobo}
\mathop{\sum_{r_1\leq x} \sum_{r_2\leq x}}_{(r_1,r_2)=1} \frac{\mu(r_1) \mu(r_2)}{
\sigma(r_1) \sigma(r_2)}.
\end{equation}
We can bound (\ref{eq:rotobo}) by the use of bounds on $\sum_{n\leq t}
\mu(n)/n$, combined with the estimation of infinite products by means
of approximations to $\zeta(s)$ for $s\to 1^+$. 
 After some additional manipulations, we
obtain a bound for $S_1(U,W)$ whose main term is at most $(3/\pi^2) (x/W)$
for each $W$, and closer to $0.22482 x/W$ on average over $W$.

(This is as good a point as any to say that, throughout, we can use
a trick in \cite{Tao} that allows us to work with odd values of integer
variables throughout, instead of letting $m$ or $n$ range over all
integers. Here, for instance, if $m$ and $n$ are restricted to
be odd, we obtain a bound of $(2/\pi^2) (x/W)$ for individual $W$,
and $0.15107 x/W$ on average over $W$. This is so even though we are losing
some cancellation in $\mu$ by the restriction.)

Let us now bound $S_2(V,W)$. This is traditionally done by Linnik's
dispersion method. However, it should be clear that the thing to do nowadays
is to use a large sieve, and, more specifically, a large sieve for primes;
such a large sieve is nothing other than a tool for estimating expressions
such as $S_2(V,W)$.
(Incidentally, even though we are trying to save every factor of $\log$ we
can, we choose not to use small sieves at all, either here or elsewhere.) 
In order to take advantage of prime support, we use Montgomery's inequality
(\cite{MR0224585}, \cite{MR0311618}; see the expositions in
\cite[pp. 27--29]{MR0337847} and  \cite[\S 7.4]{MR2061214}) combined with
Montgomery and Vaughan's large sieve with weights \cite[(1.6)]{MR0374060},
following the general procedure in \cite[(1.6)]{MR0374060}.
We obtain a bound of the form
\begin{equation}\label{eq:ursu}
\frac{\log W}{\log \frac{W}{2 q}} \left(\frac{x}{4 \phi(q)} + \frac{q W}{\phi(q)}
\right) \frac{W}{2} 
\end{equation}
on $S_2(V,W)$, where, of course, we can also choose {\em not} to gain a factor
of $\log W/2q$ if $q$ is close to or greater than $W$.

It remains to see how to gain a factor of $|\delta|$ in the major arcs,
and more specifically in $S_2(V,W)$. To explain this, let us step back and
take a look at what the large sieve is. Given a civilized 
function $f:\mathbb{Z} \to \mathbb{C}$, Plancherel's identity tells us that
\[\int_{\mathbb{R}/\mathbb{Z}} \left|\widehat{f}\left(\alpha\right)\right|^2 
d\alpha = \sum_n |f(n)|^2.\]
The large sieve can be seen as an approximate, or statistical, version of this:
for a ``sample'' of points $\alpha_1,\alpha_2,\dotsc,\alpha_k$ satisfying
$|\alpha_i-\alpha_j|\geq \beta$ for $i\ne j$, it tells us that
\begin{equation}\label{eq:rut}
\sum_{1\leq j\leq k} \left|\widehat{f}\left(\alpha_i\right)\right|^2 
\leq (X + \beta^{-1}) \sum_n |f(n)|^2,\end{equation}
 assuming that $f$ is supported on
an interval of length $X$.

Now consider $\alpha_1 = \alpha, \alpha_2 = 2\alpha, \alpha_3 = 3\alpha\dotsc
$. If $\alpha=a/q$, then the angles $\alpha_1,\dotsc,\alpha_q$ are
well-separated, i.e., they satisfy $|\alpha_i-\alpha_j|\geq 1/q$,
and so we can apply (\ref{eq:rut}) with $\beta = 1/q$. However,
$\alpha_{q+1} =\alpha_1$. Thus, if we have an outer sum of length
$L>q$ -- in (\ref{eq:katar}), we have an outer sum of length $L=x/2W$ --
we need to split it into $\lceil L/q\rceil$ blocks of length $q$,
and so the total bound given by (\ref{eq:rut}) is 
$\lceil L/q\rceil (X+q) \sum_n |f(n)|^2$. Indeed, this is what gives us
(\ref{eq:ursu}), which is fine, but we want to do better for $|\delta|$
larger than a constant.

Suppose, then, that $\alpha = a/q + \delta/x$, where $|\delta|>8$, say.
Then the angles $\alpha_1$ and $\alpha_{q+1}$ are not identical:
$|\alpha_1 - \alpha_{q+1}|\leq q |\delta|/x$. We also see that $\alpha_{q+1}$
is at a distance at least $q |\delta|/x$ from
 $\alpha_2, \alpha_3,\dotsc \alpha_q$, provided that $q |\delta|/x < 1/q$.
We can go on with $\alpha_{q+2}, \alpha_{q+3},\dotsc$, and stop only once 
there is overlap, i.e., only once we
reach $\alpha_m$ such that $m |\delta|/x \geq 1/q$. We then give all
the angles $\alpha_1,\dotsc,\alpha_m$ -- which are separated by at least
$q |\delta|/x$ from each other 
-- to the large sieve at the same time. We do this 
$\lceil L/m\rceil \leq \lceil L/(x/|\delta| q) \rceil$ times, and obtain
a total bound of $\lceil L/(x/|\delta| q) \rceil (X + x/|\delta| q)
\sum_n |f(n)|^2$, which, for $L=x/2W$, $X = W/2$, gives us about
\[\left(\frac{x}{4 Q} \frac{W}{2}
+ \frac{x}{4}\right) \log W\]
provided that $L\geq x/|\delta| q$ and, as usual, $|\alpha-a/q|\leq 1/q Q$.
This is very small compared to the trivial bound $\lesssim 
x W/8$.

What happens if $L<x/|\delta q|$? Then there is never any overlap: we
consider all angles $\alpha_i$, and give them all together to the large sieve.
The total bound is $(W^2/4 + x W/2 |\delta| q) \log W$. If $L=x/2W$
is smaller than, say, $x/3 |\delta q|$, then we see clearly that there
are non-intersecting swarms of $\alpha_i$ around the rationals $a/q$.
We can thus save a factor of $\log$ (or rather $(\phi(q)/q) \log(W/|\delta q|)$)
by applying Montgomery's inequality, which operates by strewing
displacements of the given angles (or, here, the swarms) 
around the circle to the extent possible while keeping everything 
well-separated. In this way, we obtain a bound of the form
\[\frac{\log W}{\log \frac{W}{|\delta| q}}
\left(\frac{x}{|\delta| \phi(q)} + \frac{q}{\phi(q)} \frac{W}{2}\right) \frac{W}{2}
.\]
Compare this to (\ref{eq:ursu}); we have gained a factor of $|\delta|/4$,
and so we use this estimate when $|\delta|>4$.
(In \cite{Helf}, the criterion is $|\delta|>8$, but, since there 
we have $2\alpha = a/q + \delta/x$, the value of $\delta$ there is twice what
it is here; this is a consequence of working with sums over the odd integers, 
as in \cite{Tao}.)

\begin{center}
* * *
\end{center}

We have succeeded in eliminating all factors of $\log$ we came across.
The only factor of $\log$ that remains is $\log x/UV$, coming from
the integral $\int_V^{x/U} dW/W$. Thus, we want $U V$ to be close to $x$,
but we cannot let it be too close,
since we also have a term proportional to
 $D = UV$ in (\ref{eq:marxophone}), and we need to keep it substantially
smaller than $x$. We set $U$ and $V$
so that $UV$ is $x/\sqrt{q \max(4,|\delta|)}$ or thereabouts.

In the end, after some work, we obtain the main result in \cite{Helf}.
We recall that $S_{\eta}(\alpha,x) = \sum_n \Lambda(n) e(\alpha n) \eta(n/x)$
and $\eta_2 = \eta_1 \ast_M \eta_1 = 4\cdot 1_{\lbrack 1/2,1\rbrack} \ast
1_{\lbrack 1/2,1\rbrack}$.

\begin{thm}
Let $x\geq x_0$, $x_0 = 2.16\cdot 10^{20}$.
Let $2 \alpha = a/q + \delta/x$, $q\leq Q$,
$\gcd(a,q)=1$, $|\delta/x|\leq 1/q Q$, where $Q = (3/4) x^{2/3}$.
If $q\leq x^{1/3}/6$, then
\begin{equation}\label{eq:kraw}
\begin{aligned}
&|S_{\eta}(\alpha,x)| \leq 
 \frac{R_{x,\delta_0 q} \log \delta_0 q + 0.5}{\sqrt{\delta_0 \phi(q)}} \cdot x
+ \frac{2.5 x}{\sqrt{\delta_0 q}} + 
 \frac{2x}{\delta_0 q} \cdot L_{x,\delta_0,q}
+ 3.2 x^{5/6},\end{aligned}\end{equation} where
$\delta_0 = \max(2,|\delta|/4)$,
\begin{equation}\label{eq:tosca}\begin{aligned}
R_{x,t} &= 0.27125 \log 
\left(1 + \frac{\log 4 t}{2 \log \frac{9 x^{1/3}}{2.004 t}}\right)
 + 0.41415 \\
L_{x,\delta,q} &=
\frac{\log \delta^{\frac{7}{4}} q^{\frac{13}{4}} + \frac{80}{9}}{\phi(q)/q} +
\log q^{\frac{80}{9}} \delta^{\frac{16}{9}} + \frac{111}{5}.
\end{aligned}\end{equation}
If $q > x^{1/3}/6$, then
\[|S_{\eta}(\alpha,x)|\leq 
0.2727 x^{5/6} (\log x)^{3/2} + 1218 x^{2/3} \log x.\]
\end{thm}
The factor $R_{x,t}$ is small in practice; for typical ``difficult' values
of $x$ and $\delta_0 x$, it is less than $1$. The crucial things to notice
in (\ref{eq:kraw}) are that there is no factor of $\log x$, and that,
in the main term, there is only one factor of $\log \delta_0 q$.
The fact that $\delta_0$ helps us as it grows is precisely what enables
us to take major arcs that get narrower and narrower as $q$ grows.

\section{Integrals over the major and minor arcs}\label{sec:putall}

So far, we have sketched (\S \ref{sec:rubosa}) how to estimate
$S_\eta(\alpha,x)$ for $\alpha$ in the major arcs and $\eta$ based on the
Gaussian $e^{-t^2/2}$, and also (\S \ref{sec:melusa}) how to bound $|S_\eta(\alpha,x)|$ for
$\alpha$ in the minor arcs and $\eta = \eta_2$, where $\eta_2 =
4 \cdot 1_{\lbrack 1/2,1\rbrack}\ast_M 1_{\lbrack 1/2,1\rbrack}$.
 We now must show how to use such information to 
estimate integrals such as the ones in (\ref{eq:wtree}).

We will use two smoothing functions $\eta_+$, $\eta_*$; 
in the notation of
(\ref{eq:willow}), we set $f_1 = f_2 = \Lambda(n) \eta_+(n/x)$,
 $f_3 = \Lambda(n) \eta_*(n/x)$, and so
we must give a lower bound for
\begin{equation}\label{eq:centi}
\int_{\mathfrak{M}} (S_{\eta_+}(\alpha,x))^2 S_{\eta_*}(\alpha,x) e(- \alpha n)
d\alpha\end{equation}
and an upper bound for
\begin{equation}\label{eq:wizgor}
\int_{\mathfrak{m}} \left|S_{\eta_+}(\alpha,x)\right|^2
S_{\eta_*}(\alpha,x) e(- \alpha n) d\alpha
\end{equation}
so that we can verify (\ref{eq:wtree}). 

The traditional approach to (\ref{eq:wizgor}) is to bound
\begin{equation}\label{eq:katju}\begin{aligned}
\int_{\mathfrak{m}} (S_{\eta_+}(\alpha,x))^2
S_{\eta_*}(\alpha,x) e(-\alpha n) d\alpha &\leq
\int_{\mathfrak{m}} \left|S_{\eta_+}(\alpha,x)\right|^2 d\alpha \cdot \max_{\alpha\in \mathfrak{m}}
\widehat{\eta_*}(\alpha) \\ &\leq \sum_n
\Lambda(n)^2 \eta_+^2\left(\frac{n}{x}\right) \cdot 
 \max_{\alpha\in \mathfrak{m}} S_{\eta_*}(\alpha,x).\end{aligned}\end{equation}
Since the sum over $n$ is of the order of $x \log x$, this is not
log-free, and so cannot be good enough; we will later see how to do better.
Still, this gets the main shape right: our bound on (\ref{eq:wizgor})
will be proportional to $|\eta_+|_2^2 |\eta_*|_1$. Moreover,
we see that $\eta_*$ has to be such that we know how to bound
$|S_{\eta_*}(\alpha,x)|$ for $\alpha\in \mathfrak{m}$, while our choice of
$\eta_+$ is more or less free, at least as far as the minor arcs are concerned.

What about the major arcs? In order to do anything on them, we will have
to be able to estimate both $\eta_+(\alpha)$ and $\eta_*(\alpha)$ for
$\alpha \in \mathfrak{M}$. If that is the case, then, as we shall see,
we will be able to obtain that
the main term of (\ref{eq:centi}) is
an infinite product (independent of the smoothing functions), times $x^2$,
times
\begin{equation}\label{eq:sust}\begin{aligned}
\int_{-\infty}^\infty &(\widehat{\eta_+}(- \alpha))^2 
\widehat{\eta_*}(- \alpha) e(- \alpha n/x) d\alpha \\ &=
\int_0^\infty \int_0^\infty \eta_+(t_1) \eta_+(t_2) 
\eta_*\left(\frac{n}{x} - (t_1+t_2)\right) dt_1 dt_2.
\end{aligned}\end{equation}
In other words, we want to maximize (or nearly maximize)
 the expression on the right of (\ref{eq:sust})
divided by $|\eta_+|_2^2 |\eta_*|_1$.

One way to do this is to let $\eta_*$ be concentrated on a small interval
$\lbrack 0,\epsilon)$. Then the right side of (\ref{eq:sust}) is
approximately
\begin{equation}\label{eq:korl}
\left|\eta_*\right|_1\cdot \int_0^\infty \eta_+(t) \eta_+\left(
\frac{n}{x} - t\right) dt.\end{equation}
To maximize this, we should make sure that $\eta_+(t) \sim
\eta_+(n/x-t)$. We set $x\sim n/2$, and see that we should define
$\eta_+$ so that it is supported on $\lbrack 0,2\rbrack$
and symmetric around $t=1$, or nearly so; this will maximize the ratio of
(\ref{eq:korl}) to $|\eta_+|_2^2 |\eta_*|_1$.

We should do this while making sure that we will know how to estimate
$S_{\eta_+}(\alpha,x)$ for $\alpha\in \mathfrak{M}$.
 We know how to estimate $S_\eta(\alpha,x)$ very precisely 
for functions of the form $\eta(t) = g(t) e^{-t^2/2}$, 
$\eta(t) = g(t) t e^{-t^2/2}$, etc., where $g(t)$ is band-limited. We will
work with a function $\eta_+$ of that form, chosen so as to be very
close (in $\ell_2$ norm) 
to a function $\eta_\circ$ that is in fact supported on $\lbrack 0,2\rbrack$
and symmetric around $t=1$. 

We choose
\[\eta_\circ(t) = \begin{cases}
t^2 (2-t)^3 e^{-(t-1)^2/2} & \text{if $t\in \lbrack 0,2\rbrack$,}\\
0 &\text{if $t\not\in \lbrack 0,2\rbrack$.}\end{cases}\]
This function is obviously symmetric ($\eta_\circ(t) = \eta_\circ(2-t)$)
and vanishes to high order at $t=0$, besides being supported on $\lbrack 0,2
\rbrack$.

We set $\eta_+(t) = h_R(t) t e^{-t^2/2}$, where $h_R(t)$ is an approximation to
the function \[h(t) = \begin{cases} t^2 (2-t)^3 e^{t-\frac{1}{2}}
&\text{if $t\in \lbrack 0,2\rbrack$}\\
0 &\text{if $t\not\in \lbrack 0,2\rbrack$.}\end{cases}\]
We just let $h_R(t)$ be the inverse Mellin transform of the truncation
of $Mh$ to an interval $\lbrack -iR, iR\rbrack$, or, what is the same,
\[h_R(t) = \int_0^\infty h(t y^{-1}) F_R(y) \frac{dy}{y},\]
where $F_R(t) = \sin(R \log y)/(\pi \log y)$ (the Dirichlet kernel with
a change of variables); since the Mellin transform of
$t e^{-t^2/2}$ is regular at $s=0$, the Mellin transform $M\eta_+$
will be holomorphic in a neighborhood of $\{s: 0\leq \Re(s)\leq 1\}$, even
though the truncation of $Mh$ to $\lbrack - i R, i R\rbrack$ is brutal.
Set $R=200$, say. By the fast decay of 
$Mh(it)$ and the fact that the Mellin transform $M$ is an isometry,
$|(h_R(t)-h(t))/t|_2$ is very small, and hence so is $|\eta_+-\eta_\circ|_2$,
as we desired.

But what about the requirement that we be able to
estimate $S_{\eta_*}(\alpha,x)$ for both $\alpha\in \mathfrak{m}$ and
$\alpha\in \mathfrak{M}$?

Generally speaking, if we know how to estimate
$S_{\eta_1}(\alpha,x)$ for some $\alpha\in \mathbb{R}/\mathbb{Z}$ and we
also know how to estimate $S_{\eta_2}(\alpha,x)$ for all other
$\alpha \in \mathbb{R}/\mathbb{Z}$, where $\eta_1$ and $\eta_2$ are
two smoothing functions, then we know how to estimate 
$S_{\eta_3}(\alpha,x)$ for all $\alpha\in \mathbb{R}/\mathbb{Z}$, where
$\eta_3 = \eta_1 \ast_M \eta_2$, or, more generally,
$\eta_*(t) = (\eta_1 \ast_M \eta_2)(\kappa t)$, $\kappa>0$ a constant.
This is a simple exercise in exchanging the order of integration and 
summation:
\[\begin{aligned} 
S_{\eta_*}(\alpha,x) &=
\sum_n \Lambda(n) e(\alpha n) (\eta_1\ast_M \eta_2)\left(\kappa \frac{n}{x}\right)
\\
&= \int_0^\infty \sum_n \Lambda(n) e(\alpha n) \eta_1(\kappa r) \eta_2\left(\frac{n}{r
x}\right) \frac{dr}{r} = \int_0^\infty \eta_1(\kappa r) S_{\eta_2}(r x) \frac{dr}{r},
\end{aligned}\]
and similarly with $\eta_1$ and $\eta_2$ switched. 
Of course, this trick is valid for all exponential sums:
 any function $f(n)$ would do in place of $\Lambda(n)$.
The only caveat is that
$\eta_1$ (and $\eta_2$) should be small very near $0$, since, for $r$ small,
we may not be able to estimate $S_{\eta_2}(r x)$
(or $S_{\eta_1}(r x)$) with any precision. This is not a problem; one
of our functions will be $t^2 e^{-t^2/2}$, which vanishes to second
order at $0$, and the other one will
be $\eta_2 =
4 \cdot 1_{\lbrack 1/2,1\rbrack}\ast_M 1_{\lbrack 1/2,1\rbrack}$, which has 
support bounded away from $0$. We will set $\kappa$ large (say $\kappa=49$)
so that the support of $\eta_*$ is indeed concentrated on a small interval
$\lbrack 0, \epsilon)$, as we wanted.

\begin{center}
* * *
\end{center}

Now that we have chosen our smoothing weights $\eta_+$ and $\eta_*$,
we have to estimate the major-arc integral (\ref{eq:centi}) and
the minor-arc integral (\ref{eq:wizgor}). What follows can actually be
done for general $\eta_+$ and $\eta_*$; we could have left our particular
 choice of $\eta_+$ and $\eta_*$ for the end.

Estimating the major-arc integral (\ref{eq:centi}) may
sound like an easy task, since we have rather precise estimates for
$S_{\eta}(\alpha,x)$ ($\eta=\eta_+,\eta_*$) when $\alpha$ is on the major
arcs; we could just replace $S_\eta(\alpha,x)$ in (\ref{eq:centi}) by
the approximation given by (\ref{eq:sorrow}) and (\ref{eq:frank}).
It is, however, more efficient to express (\ref{eq:centi}) as the sum of 
the contribution of the trivial character (a sum of integrals of
$(\widehat{\eta}(-\delta) x)^3$, where $\widehat{\eta}(-\delta) x$ comes
from (\ref{eq:frank})), plus a term of the form
\[(\text{maximum of $\sqrt{q}\cdot E(q)$ for $q\leq r$})\cdot
\int_{\mathfrak{M}} \left|S_{\eta_+}(\alpha,x)\right|^2 d\alpha
,\]
where $E(q)=E$ is as in (\ref{eq:lynno}), plus two other terms of the same form.
As usual, the major arcs $\mathfrak{M}$ are the arcs around rationals
$a/q$ with $q\leq r$.
We will soon discuss how to bound the integral of 
$\left|S_{\eta_+}(\alpha,x)\right|^2$ over arcs around rationals $a/q$ with
 $q\leq s$,
$s$ arbitrary. Here, however, it is best to estimate the integral over 
$\mathfrak{M}$ using the estimate on
$S_{\eta_+}(\alpha,x)$ from (\ref{eq:sorrow}) and (\ref{eq:frank}); we obtain
a great deal of cancellation, with the effect that,
for $\chi$ non-trivial, the error term in 
(\ref{eq:lynno}) appears only when it gets squared,
and thus becomes negligible.

The contribution of the trivial character has an easy approximation, 
thanks to the fast decay of $\widehat{\eta_\circ}$.
We obtain that the major-arc integral (\ref{eq:centi}) equals a main
term $C_0 C_{\eta_\circ,\eta_*} x^2$, where
\[\begin{aligned}
C_0 &= \prod_{p|n} \left(1 - \frac{1}{(p-1)^2}\right)
\cdot \prod_{p\nmid n} \left(1 + \frac{1}{(p-1)^3}\right),\\
C_{\eta_\circ,\eta_*} &=
 \int_0^\infty \int_0^\infty \eta_\circ(t_1) \eta_\circ(t_2) 
\eta_*\left(\frac{n}{x}-(t_1+t_2)\right) dt_1 dt_2,
\end{aligned}\]
plus several small error terms. We have already chosen $\eta_\circ$,
$\eta_*$ and $x$ so as
to (nearly) maximize $C_{\eta_\circ,\eta_*}$.

It is time to bound the minor-arc integral (\ref{eq:wizgor}).
As we said in \S \ref{sec:putall}, we must do better than the usual bound
(\ref{eq:katju}). Since our minor-arc bound (\ref{eq:kraw})
on $|S_\eta(\alpha,x)|$, $\alpha\sim a/q$, decreases as $q$ increases,
it makes sense to use partial summation together with bounds on
\[\int_{\mathfrak{m}_s}  |S_{\eta_+}(\alpha,x)|^2
= \int_{\mathfrak{M}_s} |S_{\eta_+}(\alpha,x)|^2
d\alpha -
\int_{\mathfrak{M}} |S_{\eta_+}(\alpha,x)|^2
d\alpha,
\]
where
$\mathfrak{m}_s$ denotes the arcs around $a/q$, $r<q\leq s$, and
$\mathfrak{M}_s$ denotes the arcs around all $a/q$, $q\leq s$.  We already
know how to estimate the integral on $\mathfrak{M}$. How do we bound the
integral on $\mathfrak{M}_s$?

In order to do better than the trivial bound $\int_{\mathfrak{M}_s} \leq
\int_{\mathbb{R}/\mathbb{Z}}$,
we will need to use the fact that the series (\ref{eq:holo}) 
defining $S_{\eta_+}(\alpha,x)$ is essentially supported on prime numbers.
Bounding the integral on $\mathfrak{M}_s$ is closely related to the problem
of bounding
\begin{equation}\label{eq:kolro}
\sum_{q\leq s} \mathop{\sum_{a \mo q}}_{(a,q)=1} \left|
\sum_{n\leq x} a_n e(a/q)\right|^2
\end{equation}
efficiently for $s$ considerably smaller than $\sqrt{x}$ and
$a_n$ supported on the primes $\sqrt{x}<p\leq x$. This is a classical
problem in the study of the large sieve. The usual bound on (\ref{eq:kolro})
(by, for instance, Montgomery's inequality) has a gain of a factor of
$2 e^{\gamma} (\log s)/(\log x/s^2)$ relative to the bound of $(x+s^2)
\sum_n |a_n|^2$ that one would get from the large sieve without using prime 
support. Heath-Brown proceeded similarly to bound
\begin{equation}\label{eq:kokoto}\int_{\mathfrak{M}_s} |S_{\eta_+}(\alpha,x)|^2
d\alpha 
\lesssim \frac{2 e^{\gamma} \log s}{\log x/s^2}
\int_{\mathbb{R}/\mathbb{Z}} |S_{\eta_+}(\alpha,x)|^2
d\alpha.\end{equation}

This already gives us the gain of
$C (\log s)/\log x$ that we absolutely need,
but the constant $C$ is suboptimal; the factor in the right side of
(\ref{eq:kokoto}) should really be
$(\log s)/\log x$, i.e., $C$ should be $1$. We cannot reasonably hope to do better than $2 (\log s)/\log x$ in the minor arcs due to what is known as the {\em parity problem} 
in sieve theory. As it turns out, Ramar\'e \cite{MR2493924} had given general
bounds on the large sieve that were clearly conducive to better bounds
on (\ref{eq:kolro}), though they involved a ratio that was not easy to
bound in general.

I used several careful estimations (including \cite[Lem.~3.4]{MR1375315})
to reduce the problem of bounding this ratio to a finite number of cases, which
 I then checked by rigorous computation. This approach gave a bound on
(\ref{eq:kolro}) with a factor of size close to $2 (\log s)/\log x$.
(This solves the large-sieve problem for $s\leq x^{0.3}$; it would still
be worthwhile to give a computation-free proof for all $s\leq x^{1/2-\epsilon}$,
$\epsilon>0$.) It was then easy to give an analogous bound for the 
integral over $\mathfrak{M}_s$, namely,
\[\int_{\mathfrak{M}_s} |S_{\eta_+}(\alpha,x)|^2
d\alpha 
\lesssim \frac{2 \log s}{\log x}
\int_{\mathbb{R}/\mathbb{Z}} |S_{\eta_+}(\alpha,x)|^2
d\alpha,\]
where $\lesssim$ can easily be made precise by replacing $\log s$ by
$\log s + 1.36$ and $\log x$ by $\log x + c$, where $c$ is a small constant.
Without this improvement, the main theorem would still have been proved, but
the required computation time would have been multiplied by a factor of
considerably more than $e^{3\gamma} = 5.6499\dotsc$. 

What remained then was just to compare the estimates on 
(\ref{eq:centi}) and (\ref{eq:wizgor}) and check that (\ref{eq:wizgor})
is smaller for $n\geq 10^{27}$. This final step was just bookkeeping.
As we already discussed, a check for $n< 10^{27}$ is easy.
Thus ends the proof of the main theorem.
\section{Some remarks on computations}

There were two main computational tasks: verifying the ternary conjecture
for all $n\leq C$, and checking the Generalized Riemann Hypothesis
for modulus $q\leq r$ up to a certain height. 

The first task was not very
demanding. Platt and I verified in \cite{HelPlat}
 that every odd integer 
$5 < n\leq 8.8\cdot 10^{30}$
can be written as the sum of three primes. (In the end, only a check
for $5<n\leq 10^{27}$ was needed.) 
We proceeded as follows. In a major computational effort,
Oliveira e Silva, Herzog and Pardi \cite{OSHP}) had already checked
that the binary Goldbach conjecture is true up to $4\cdot 10^{18}$ --  
that is, every even number up to $4\cdot 10^{18}$ is the sum of two primes. 
Given that, all we had to do was to construct a ``prime ladder'', that is, a 
list of primes from $3$ up to $8.8\cdot 10^{30}$ such that the difference between any two consecutive primes in the list is at least $4$ and
at most $4\cdot 10^{18}$. 
(This is a known strategy: see \cite{MR1451327}.) Then,
for any odd integer $5 < n\leq 8.8\cdot 10^{30}$, there is a prime $p$ in 
the list such that $4\leq n-p \leq 4\cdot 10^{18}+2$. (Choose the largest $p<n$
in the ladder, or, if $n$ minus that prime is $2$, choose the prime immediately
under that.)
By \cite{OSHP} (and the fact that $4\cdot 10^{18}+2$ equals $p+q$, where
$p=2000000000000001301$ and $q=1999999999999998701$ are both prime), 
we can write $n-p = p_1 + p_2$ for some primes $p_1$, $p_2$, 
and so $n = p + p_1 + p_2$.

Building a prime ladder involves only integer arithmetic, that is, computer
manipulation of integers, rather than of real numbers.
Integers are something that computers can handle rapidly and reliably.
We look for primes for our ladder only among a special set of integers 
whose primality can be tested deterministically quite quickly (Proth numbers: 
$k\cdot 2^m+1$, $k< 2^m$). Thus, we can build a prime ladder by
 a rigorous, deterministic
algorithm that can be (and was) parallelized trivially.

The second computation is more demanding. It consists in 
verifying that, for every $L$-function $L(s,\chi)$ with $\chi$ 
of conductor $q\leq r = 300000$ (for $q$ even) or $q\leq r/2$ (for $q$ odd), 
all zeroes of $L(s,\chi)$ such that $|\Im(s)|\leq H_q = 10^8/q$ (for $q$ odd)
and $|\Im(s)|\leq H_q = \max(10^8/q,200+7.5\cdot 10^7/q$ (for $q$ even)
 lie on the critical line. 
This was entirely Platt's work; my sole contribution was to request computer
time.
In fact, he went up to conductor $q\leq 200000$ (or twice that for $q$ even); 
he had already gone up to conductor $100000$ in his PhD thesis. The 
verification took, in total, about $400000$ core-hours (i.e., the total number 
of processor cores used times the number of hours they ran equals $400000$; 
nowadays, a top-of-the-line processor typically has eight cores). In the end, 
since I used only $q\leq 150000$ (or twice that for $q$ even), 
the number of hours actually needed was closer to $160000$; since 
I could have made do with $q\leq 120000$ (at the cost of increasing $C$ to $10^{29}$
or $10^{30}$), it is likely, in retrospect, that only about $80000$ 
core-hours were needed.

Checking zeros of $L$-functions computationally goes back to Riemann 
(who did it by hand for the special case of the Riemann zeta function). 
It is also one of the things that were tried on digital computers in 
their early days (by Turing \cite{MR0055785}, for instance; see the 
exposition in \cite{MR2263990}). One of the main issues to be careful about 
arises whenever one manipulates real numbers via a computer: generally speaking,
a computer cannot store an irrational number; moreover, while a computer can handle rationals, it is really most comfortable handling just those rationals whose denominators are powers of two. Thus, one cannot really say: ``computer, give me the sine of that number'' and expect a precise result. What one should do, if one really wants to prove something (as is the case here!), is to say: 
``computer, I am giving you an interval $I=\lbrack a/2^k,b/2^k\rbrack$; 
give me an interval $I'=\lbrack c/2^\ell,d/2^\ell
\rbrack$, preferably very short, such that $\sin(I) \subset I'$''. This is called interval arithmetic; it is arguably the easiest way to do floating-point computations rigorously.

Processors do not do this natively, and if interval arithmetic is
implemented purely on software, computations can be slowed down by a factor of
about $100$. Fortunately, there are ways of running interval-arithmetic 
computations
partly on hardware, partly on software. Platt has his own library, but there 
are others online (e.g. PROFIL/BIAS \cite{Profbis}).

Incidentally, there are some basic functions (such as $\sin$) that should
always be done on software, not just if one wants to use interval arithmetic, but
even if one just wants reasonably precise results: the implementation of
transcendental functions in some of the most popular processors (Intel)
does not always round correctly, and errors can accumulate quickly.
Fortunately, this problem is already well-known, and there is software
(for instance, the crlibm library \cite{crlibm}) that takes care of this.

Lastly, there were several relatively minor computations embedded in
\cite{Helf}, \cite{HelfMaj}, \cite{HelfTern}. There is some numerical 
integration, done rigorously; this is sometimes done using a standard
package based on interval arithmetic \cite{VNODELP}, but most of the time
I wrote my own routines in C (using Platt's interval arithmetic package) for
the sake of speed. Another typical computation was a 
rigorous version of a ``proof by graph'' (``the maximum of a function $f$ is clearly less than $4$ because I can see it on the screen''). There is a 
standard way to do this (see, e.g.,  \cite[\S 5.2]{MR2807595}); 
essentially, the bisection method combines naturally with interval arithmetic.
Yet another computation (and not a very small one)
was that involved in verifying a large-sieve 
inequality in an intermediate range (as we discussed in \S \ref{sec:putall}).

 It may be interesting to note that one of the 
inequalities used to estimate (\ref{eq:rotobo})
was proven with the help of automatic quantifier elimination
\cite{QEPCAD}. Proving this inequality was a very 
minor task, both computationally and mathematically;
 in all likelihood, it is feasible to give a human-generated 
proof. Still, it is nice to know from first-hand experience
that computers can nowadays (pretend to) do something
other than just perform numerical computations -- and that this is true
even in current mathematical practice.
\bibliographystyle{alpha}
\bibliography{icm}
\end{document}